\documentclass[leqno, 12pt]{amsart}
\usepackage{amscd}
\usepackage{amssymb}
\usepackage[all, poly]{xy}
\newcommand{\nc}{\newcommand}
\newcommand{\rc}{\renewcommand}

\nc{\ra}{{      \rightarrow     }}
\nc{\laa}{{     \leftarrow      }}      %Forget it and use @<<<
\nc{\lra}{{\longrightarrow}}
\nc{\lr}{{\leftrightarrow}}             % "corresponds"
\nc{\lrs}{{\rightleftarrows}}           % One map in each direction.

\nc{\imp}{{\Rightarrow}}                %Implication.
\nc{\eq}{{\Leftrightarrow}}             %Equivalence.

\nc{\inj}{{\pr  \hookrightarrow }}              %Injective map-right.
\nc{\injj}{{\pr \hookleftarrow  }}              %Injective map.

\nc{\sur}{{     \twoheadrightarrow      }}      %Surjective map-right.
\nc{\surr}{{    \twoheadleftarrow       }}      %Surjective map-left

\nc{\mm}{{\mapsto}}                             %Map on elements.

\nc{\va}{{\uparrow}}                            %Vertical Up-arrow
\nc{\barr}{\overline}
\nc{\ul}{\underline}
\nc{\sub}{\subseteq}

\nc{\se}{       \section                }
\nc{\sus}{      \subsection             }
\nc{\sss}{      \subsubsection          }

\nc{\Lemm}{     \subsection{Lemma}              }
\nc{\lemm}{     \subsubsection{Lemma}           }
\nc{\slemm}{    \subsubsection*{Lemma}          }
\nc{\sublemm}{  \subsubsection{\bf Sublemma}            }
\nc{\ssublemm}{         \subsubsection*{\bf Sublemma}           }

\nc{\Pro}{      \subsection{Proposition}        }
\nc{\pro}{      \subsubsection{Proposition}     }
\nc{\spro}{     \subsubsection*{Proposition}    }

\nc{\Corr}{     \subsection{Corollary}          }
\nc{\corr}{     \subsubsection{\bf Corollary}       }
\nc{\scorr}{    \subsubsection*{\bf Corollary}      }

\nc{\Theo}{     \subsection{Theorem}            }
\nc{\theo}{     \subsubsection{\rm{\bf Theorem}} }
\nc{\stheo}{    \subsubsection*{\rm{\bf Theorem}}        }

\nc{\rem}{      \subsubsection{Remark}          }
\nc{\srem}{     \subsubsection*{Remark} }

\nc{\rems}{     \subsubsection{Remarks}         }
\nc{\srems}{    \subsubsection*{Remarks}        }
\nc{\conj}{     \subsubsection{\bf Conjecture}      }
\nc{\sconj}{    \subsubsection*{\bf Conjecture}     }

\nc{\ex}{       \subsubsection{Example}         }
\nc{\sex}{      \subsubsection*{Example}        }
\nc{\exs}{      \subsubsection{Examples}        }
\nc{\sexs}{     \subsubsection*{Examples}       }

\nc{\que}{      \subsubsection{Question}        }
\nc{\ques}{     \subsubsection{Questions}       }
\nc{\sque}{     \subsubsection*{Question}       }
\nc{\sques}{    \subsubsection*{Questions}      }

\nc{\pl}{{\oplus}}                              %Direct sum
\nc{\tim}{{\times}}
\nc{\btim}{{\boxtimes}}
\nc{\ltim}{\ltimes}                     %
\nc{\rtim}{\rtimes}                     %
\nc{\ltr}{\triangleleft}        %
\nc{\rtr}{\triangleright}       %

\nc{\ten}{{     \otimes         }}
\nc{\Lten}{{    \aa{L}\otimes   }}            %Derived tensoring.
\nc{\Ltim}{{    \aa{L}\times    }}            %Derived product.
\nc{\bten}{{\boxtimes}}                         %Tensoring: outer

\nc{\con}{{ @>{\protect\cong}>> }}      %Isomorphism with a right arrow
\nc{\conl}{{    @>{\cong}>>     }}      %Lower right isomorphism
\nc{\conn}{{    @<{\cong}<<     }}      %Isomormphism with a left arrow
\nc{\Con}{{     \equiv          }}      %Congruence
\nc{\appr}{{    \sim            }}      %Approximately
\nc{\eqr}{{     \sim            }}      %Equivalence relation

\nc{\ha}{{ \frac{1}{2} }}               %%%Half
        \nc{\half}{{ \frac{1}{2} }}

\nc{\ci}{{\circ}}               %Circle dot
\nc{\cd }{{\cdot}}              %Multiplication dot
\nc{\cddd}{{\cdots}}

\nc{\cupp}{\bigcup}             %cup,cap
\nc{\capp}{\bigcap}
\nc{\pll}{\bigoplus}

\nc{\pii}{\prod}                %Product
\nc{\ppii}{\bigprod}            %Big product
\nc{\cci}{\sqcup}              %Coprodduct
\nc{\ccii}{\bigsqcup}

\nc{\wwe}{\bigwedge}            %wedge
\nc{\cce}{\bigcoprod}           %cowedge

\nc{\aaa}{      \stackerel      }       %Stack {1} over 2.

\rc{\Im}{{      \operatorname{Im}       }}
\nc{\rank}{{    \ \operatorname{rank}\  }}
\nc{\Res}{{     \  \operatorname{Res}   }}
\nc{\Hom}{{    \operatorname{Hom}      }}
\nc{\End}{{     \operatorname{End}      }}
\nc{\RHom}{{    \operatorname{RHom}     }}
\nc{\HHom}{{    \operatorname{$\HH$om}  }}
\nc{\EEnd}{{    \operatorname{$\EE nd$} }}
\nc{\AAut}{{    \operatorname{$\AA ut$} }}
\nc{\RHHom}{{   \operatorname{R$\HH$om} }}
\nc{\Ext}{\operatorname{Ext}}
\nc{\Der}{{     \operatorname{Der}      }}

\nc{\ord        }{{ \operatorname{ord} }}                       %order of zero
\nc{\divv       }{{ \operatorname{div} }}                       %divisor
\nc{\Lie        }{{ \operatorname{Lie} }}

\rc{\AA}{{\mathcal A}}
\nc{\BB}{{\mathcal B}}
\nc{\CC}{{\mathcal C}}
\nc{\DD}{{\mathcal D}}
\nc{\EE}{{\mathcal E}}
\nc{\FF}{{\mathcal F}}
\nc{\GG}{{\mathcal G}}
\nc{\HH}{{\mathcal H}}
\nc{\II}{{\mathcal I}}
\nc{\JJ}{{\mathcal J}}
\nc{\KK}{{\mathcal K}}
\nc{\LL}{{\mathcal L}}
\nc{\MM}{{\mathcal M}}
\nc{\NN}{{\mathcal N}}
\nc{\OO}{{\mathcal O}}
\nc{\PP}{{\mathcal P}}
\nc{\QQ}{{\mathcal Q}}
\nc{\RR}{{\mathcal R}}
\rc{\SS}{{\mathcal S}}
\nc{\UU}{{\mathcal U}}
\nc{\VV}{{\mathcal V}}
\nc{\WW}{{\mathcal W}}
\nc{\ZZ}{{\mathcal Z}}
\nc{\XX}{{\mathcal X}}
\nc{\YY}{{\mathcal Y}}

\nc{\A}{{\mathbb A }}
\nc{\B}{{\mathbb B}}
\nc{\C}{{\mathbb C}}
\nc{\D}{{\mathbb D}}
\nc{\E}{{\mathbb E}}
\nc{\F}{{\mathbb F}}
\nc{\G}{{\mathbb G}}
\nc{\hH}{{\mathbb H}}
\nc{\I}{{\mathbb I}}
\nc{\J}{{\mathbb J}}
\nc{\K}{{\mathbb K}}
\nc{\lL}{{\mathbb L}}
\nc{\M}{{\mathbb M}}
\nc{\N}{{\mathbb N}}
\nc{\oO}{{\mathbb O}}
\nc{\pP}{{\mathbb P}}
\nc{\Q}{{\mathbb Q}}
\nc{\R}{{\mathbb R}}
\nc{\sS}{{\mathbb S}}
\nc{\T}{{\mathbb T}}
\nc{\U}{{\mathbb U}}
\nc{\V}{{\mathbb V}}
\nc{\W}{{\mathbb W}}
\nc{\Z}{{\mathbb Z}}
\nc{\X}{{\mathbb X}}
\nc{\Y}{{\mathbb Y}}

\nc{\fA}{{\mathfrak A}}
\nc{\fB}{{\mathfrak B}}
\nc{\fC}{{\mathfrak C}}
\nc{\fD}{{\mathfrak D}}
\nc{\fE}{{\mathfrak E}}
\nc{\fF}{{\mathfrak F}}
\nc{\fG}{{\mathfrak G}}
\nc{\fH}{{\mathfrak H}}
\nc{\fI}{{\mathfrak I}}
\nc{\fJ}{{\mathfrak J}}
\nc{\fK}{{\mathfrak K}}
\nc{\fL}{{\mathfrak L}}
\nc{\fM}{{\mathfrak M}}
\nc{\fN}{{\mathfrak N}}
\nc{\fO}{{\mathfrak O}}
\nc{\fP}{{\mathfrak P}}
\nc{\fQ}{{\mathfrak Q}}
\nc{\fR}{{\mathfrak R}}
\nc{\fS}{{\mathfrak S}}
\nc{\fT}{{\mathfrak T}}
\nc{\fU}{{\mathfrak U}}
\nc{\fV}{{\mathfrak V}}
\nc{\fW}{{\mathfrak W}}
\nc{\fZ}{{\mathfrak Z}}
\nc{\fX}{{\mathfrak X}}
\nc{\fY}{{\mathfrak Y}}
\nc{\fa}{{\mathfrak a}}
\nc{\fb}{{\mathfrak b}}
\nc{\fc}{{\mathfrak c}}
\nc{\fd}{{\mathfrak d}}
\nc{\fe}{{\mathfrak e}}
\nc{\ff}{{\mathfrak f}}
\nc{\fg}{{\mathfrak g}}
\nc{\fh}{{\mathfrak h}}
\nc{\fj}{{\mathfrak j}}
\nc{\fk}{{\mathfrak k}}
\nc{\fl}{{\mathfrak{l}}}
\nc{\fm}{{\mathfrak m}}
\nc{\fn}{{\mathfrak n}}
\nc{\fo}{{\mathfrak o}}
\nc{\fp}{{\mathfrak p}}
\nc{\fq}{{\mathfrak q}}
\nc{\fr}{{\mathfrak r}}
\nc{\fs}{{\mathfrak s}}
\nc{\ft}{{\mathfrak t}}
\nc{\fu}{{\mathfrak u}}
\nc{\fv}{{\mathfrak v}}
\nc{\fw}{{\mathfrak w}}
\nc{\fz}{{\mathfrak z}}
\nc{\fx}{{\mathfrak x}}
\nc{\fy}{{\mathfrak y}}

\nc{\al}{{\alpha }}
\nc{\be}{{\beta }}
\nc{\ga}{{\gamma }}
\nc{\de}{{\delta }}
\nc{\del}{{\partial }}
\nc{\ep}{{\varepsilon }}
\nc{\vap}{{\epsilon }}

\nc{\ze}{{\zeta }}
\nc{\et}{{\eta }}
\rc{\th}{{\theta }}
\nc{\vth}{{\vartheta }}

\nc{\io}{{\iota }}
\nc{\ka}{{\kappa }}
\nc{\la}{{\lambda }}
\nc{\vrho}{{\varrho}}
\nc{\si}{{\sigma }}
\nc{\ups}{{\upsilon }}
\nc{\vphi}{{\varphi }}
\nc{\om}{{\omega }}

\nc{\Ga}{{\Gamma }}
\nc{\De}{{\Delta }}
\nc{\nab}{{\nabla}}
\nc{\Th}{{\Theta }}
\nc{\La}{{\Lambda }}
\nc{\Si}{{\Sigma }}
\nc{\Ups}{{\Upsilon }}
\nc{\Om}{{\Omega }}

\nc{\toc}{{\tableofcontents}}
\nc{\addl}{     \addcontentsline{toc}{subsection}       }

\begin{document}

                        %%GRASSMANNIAN SYMBOLS
% GROUPS
                        %G:
\nc{\GK}{{      G(\KK)          }}
\nc{\GO}{{      G(\OO)          }}

\nc{\Gd}{{  {\check G}          }}

\nc{\Gh}{{  \hat{\GG}           }}
\nc{\GA}{{  G(\AA)              }}
                        %B:
\nc{\BK}{{  B(\mathcal K)           }}
\nc{\BO}{{  B(\mathcal O)           }}
\nc{\BKz}{{  B(\mathcal K)_0        }}
                        %N:
\nc{\NK}{{  N(\mathcal K)           }}
\nc{\NO}{{  N(\mathcal O)           }}
                        %T
\nc{\TK}{{  B(\mathcal K)           }}
\nc{\TO}{{  T(\mathcal O)           }}
                        %L
\nc{\LO}{{  L(\mathcal O)           }}

%LIE ALGEBRAS
\nc{\gk}{{      \mathfrak g_\KK             }}
\nc{\bk}{{      \mathfrak b_\KK             }}
\nc{\tk}{{      \mathfrak b_\KK             }}
\nc{\nk}{{      \mathfrak n_\KK             }}
\nc{\go}{{      \mathfrak g_\OO             }}
%\nc{\bo}{{     \mathfrak b_\OO             }}
%\nc{\ga}{{     \fg_\mathcal A              }}

%POINTS in Grassmannian
\nc{\Ll}{{  L_\lambda   }}
\nc{\Lm}{{  L_\mu       }}

%GO-orbits
\nc{\Gl}{{      {\mathcal G_\lambda}        }}
\nc{\Glb}{{     \barr{\Gl}              }}
\nc{\Glm}{{     \GG_{\lambda+\mu}       }}
\nc{\Ge}{{      \GG_\eta                }}
\nc{\Geb}{{     \barr{\Ge}              }}
\nc{\Gwl}{{     \GG_{W \lambda}         }}
\nc{\Gn}{{      \GG_\nu                 }}
\nc{\Gnb}{      \barr{\GG_\nu           }}
\nc{\Gm}{{      \GG_\mu                 }}
                                %Upper orbits
\nc{\Gum}{{     \GG^\mu         }}
\nc{\Gul}{{     \GG^\lambda     }}
\nc{\Gun}{{     \GG^\nu         }}
                                %STRUCTURE OF GO-ORBITS
\nc{\Pl}{{  \PP_\lambda         }}
\nc{\Pwl}{{  {\PP}_{W \lambda}  }}
\nc{\FFl}{{  {\FF}_\lambda      }}

%IWAHORI ORBITS
\nc{\Cn}{{  C_\nu       }}
\nc{\Ce}{{  C_\eta      }}

%N_K ORBITS
\nc{\Sn }{{     S_\nu           }}
\nc{\Snb}{{  \barr{S_\nu}       }}
\nc{\Sl }{{  S_\lambda          }}
\nc{\Sm }{{  S_\mu              }}

%DOUBLE Grassmannian
\nc{\GGG}{{  G(\KK)\underset{G(\OO)}\times      \GG                     }}
\nc{\GGlm}{{  (G(\KK)\underset{G(\OO)}\times    \GG)_{(\lambda,\mu)}    }}

\nc{\dgg}{{     \ddot\GG                }}
\nc{\dlm}{{     \ddot\GG_{\la,\mu}      }}

%ROOTS
\nc{\db}{{      \hat \Delta     }}
\nc{\dr}{{      \Delta_\R       }}
\nc{\Dp}{{      \Delta^+        }}
\nc{\Dm}{{      \Delta^-        }}
\nc{\da}{{      \Delta_a        }}
\nc{\dap}{{     \Delta_a^+      }}
\nc{\di}{{      \Delta_I        }}
\nc{\cro}{{     \check \rho     }}      %Former         "\rc"
\nc{\dc}{{      \check \delta   }}
%\nc{\raa}{{    \rho_a          }}           %Used to be ``\ra''
\nc{\rac}{{     \check \rho_a   }}

%%%%%%%%%%%%%%%%%%%%%%%%%%%%%%%% MAX %%%%%%%%%%%%%%%%%%%%%%%%%%%%%%%

\nc{\Irr}{\operatorname {Irr}}
\nc{\td}{\widetilde d}
\nc{\tg}{\widetilde g}
\nc{\tv}{\widetilde v}
\nc{\tp}{\tilde p}
\nc{\tmm}{\widetilde{\bf m}}
\nc{\tmu}{\widetilde{\mu}}

\nc{\TD}{\widetilde D}
\nc{\TV}{\widetilde V}
\nc{\TU}{\widetilde U}
\nc{\TT}{\widetilde T}
\nc{\TA}{\widetilde A}
\nc{\TB}{\widetilde B}
\nc{\TN}{\widetilde{\mathcal N}}
\nc{\TF}{\widetilde{\mathcal F}}
\nc{\TGG}{\widetilde{\mathcal G}}

\nc{\tga}{\widetilde\ga}
\nc{\tde}{\widetilde\de}
\nc{\tphi}{\widetilde\phi}
\nc{\tpsi}{\widetilde\psi}

\nc{\cla}{\check\lambda}
\nc{\cmu}{\check\mu}
\nc{\cth}{\check\theta}

\nc{\mmm}{{\bf m}}

\nc{\bx}{\barr x}

\nc{\Sym}{\operatorname {Sym}}
\nc{\Spec}{\operatorname {Spec}}
\nc{\Id}{\operatorname {Id}}

\nc{\gl}{\mathfrak{gl}}
\nc{\tfg}{\widetilde{\mathfrak g}}

\nc{\TBDG}{\widetilde{\bf{\mathfrak G}}}
\nc{\BDG}{{\bf\mathfrak G}}

\nc{\Perv}{\operatorname {Perv}}
\nc{\Rep}{\operatorname {Rep}}
\nc{\Mlt}{\operatorname {Mlt}}
\nc{\spr}{\operatorname {Spr}}

\nc{\IC}{\operatorname {IC}}
\nc{\Tr}{\operatorname {Tr}}
\nc{\sign}{\operatorname{sgn}}

\nc{\ad}{\operatorname{ad}}
\nc{\pr}{\operatorname{pr}}

\nc{\bfb}{{\bf b}}
\nc{\bfg}{{\bf g}}
\nc{\bfe}{{\bf e}}

%%%%%%%%%%%%%%%%%%%%%%%%% END MAX SYMBOLS %%%%%%%%%%%%%%%%%%%%%%%

                        %TITLE and AUTHORS
\title[]{
Quiver varieties and Beilinson-Drinfeld Grassmannians of type A
}

\author{       Ivan Mirkovi\'c                 }
\address{Dept. of Mathematics and Statistics, University
of Massachusetts at Amherst, Amherst MA 01003-4515, USA}
\email{                mirkovic@math.umass.edu        }

\author{       Maxim Vybornov                }
\address{ 22 Rockland St, Newton, MA 02458}
\email{                mail@maximvybornov.net }

\begin{abstract}
We construct Nakajima's quiver varieties of type A in terms of
conjugacy classes of matrices and (non-Slodowy's) transverse slices
naturally arising from affine Grassmannians. In full generality quiver varieties
are embedded into Beilinson-Drinfeld Grassmannians of type A.
Our construction provides a compactification of Nakajima's
quiver varieties and a decomposition of an affine Grassmannian
into a disjoint union of quiver varieties. As an application we
provide a geometric version of skew and symmetric $(GL(m), GL(n))$
duality.
\end{abstract}

%\thanks{}

%\date{
%Dec 24, 2007 }

\maketitle
                        %%%%%%Home%%%%
%Local definitions

%\toc

\se{Introduction}

In type A we relate Nakajima's quiver varieties, conjugacy classes of matrices,
and Beilinson-Drinfeld Grassmannians. In
particular, we embed quiver varieties into Beilinson-Drinfeld
Grassmannians.
From the point of view of Nakajima's quiver varieties our construction
provides a compactification of quiver varieties.
From the point of view of nilpotent orbits we construct new transverse slices
to nilpotent orbits naturally arising from affine Grassmannians.
From the point of view of affine Grassmannians we get a decomposition of an affine Grassmannian
into a disjoint union of quiver varieties.
As an application we provide a geometric version of both the skew and the symmetric
version of the $(GL(m), GL(n))$ duality.

The relationship between quiver varieties and
nilpotent orbits was conjectured by Nakajima \cite{N94} and
proved by Maffei \cite{M}. What we do here is close to (and in part motivated
by) Maffei's work, however while he uses
Slodowy's normal slices to nilpotent orbits we use different
slices suggested by the relation to the affine Grassmannians, and
this makes the construction explicit while Maffei's approach is
based on an existence result.

These observations do  not literally extended
beyond type A. For instance, the closures of orbits in the affine Grassmannian
are normal and this is not true for the nilpotent orbits.

\sus{The setup} We work over the field of complex numbers $\C$.
By $G_m=\C^{\ast}$ we sometimes denote the multiplicative group of this field.

Given two $(n-1)$-tuples of integers
$d=(d_1,\dots ,d_{n-1})$ and
$v=(v_1,\dots , v_{n-1})$
and a central element $c=(c_1,\dots, c_{n-1})$
of the Lie algebra $\prod_{i=1}^{n-1} \gl(v_i, \C)$,
Nakajima \cite{N94, N98} constructs quiver varieties $\fM_0(v,d)$ and
$\fM(v,d)$.

From the quiver data one can produce $GL(m)$-(co)weights (partitions)
$\la$ and $\mu$ of $N$ (cf. subsection \ref{qdatatogl}),
where $m=d_1+\dots +d_{n-1}$, and $N = \sum_{j = 1}^{n-1} jd_j$.
We will also consider the affine Grassmannian $\GG$
associated to the group $G=GL(m)$, and a ``convolution''
Grassmannian $\widetilde\GG$ equipped with a resolution map $\pi: \widetilde\GG\to \GG$.

The following theorem is a common generalization of
(some of) the results of Kraft-Procesi \cite{KP}, Lusztig \cite{L81},
and Nakajima \cite{N94}. For simplicity we will only
write down here the statement in the case $c=0$.
In this paper we provide a complete proof of the Theorem below announced in
\cite{MVyb}.

\subsection{\bf Theorem}
There exist algebraic isomorphisms
$\phi, \widetilde\phi, \psi, \widetilde\psi$
such that the following diagram commutes:
\begin{equation}\nonumber
\begin{CD}
\fM(v,d)
@>{\widetilde\phi}>{\simeq}>
{\bf m}^{-1}(T_{\la}\cap\barr\OO_\mu)
@>{\widetilde\psi}>{\simeq}>
\pi^{-1}(L^{<0}G\cdot \la\cap\barr{L^{\geq 0}G\cdot {\mu}})
&\  \subset \ & \widetilde\GG
\\
@V{p}VV
@V{\bf m}VV
@V{\pi}VV
\\
\fM_0(v,d)
@>{\phi}>{\simeq}>
T_{\la}\cap\barr\OO_\mu
@>{\psi}>{\simeq}>
L^{<0}G\cdot \la\cap\barr{L^{\geq 0}G\cdot {\mu}}
&\  \subset \ & \GG ,
\end{CD}
\end{equation}
where $T_\la$ is our new transverse slice to the nilpotent orbit $\OO_\la \sub \NN$
of type $\la$ in the nilpotent cone $\NN$ of the $\gl(N, \C)$, $\barr\OO_\mu$ is the closure
of the nilpotent orbit of type $\mu$ in $\NN$,
${\bf m}: \widetilde\OO_\mu\to \barr\OO_\mu$
is its Springer resolution, and
$L^{\geq 0}G$ and $L^{<0}G$ are the subgroups of non-negative
and negative loops respectively in the loop group $GL(m,\C((z)))$.

\sus{The deformation} For arbitrary $c$, the nilpotent orbits deform to general
conjugacy classes, and the affine Grassmannian
deforms to the Beilinson-Drinfeld Grassmannian
$\BDG_{{\mathbb A}^{(n)}}$
on the $n$-th symmetric power of the curve ${\mathbb A}^1$,
or more precisely its fiber
over the point
$(0,c_1,c_1+c_2,\dots,c_1+\dots+c_{n-1})\in {\mathbb A}^{(n)}$.
The general statement is formulated as Theorem \ref{main}.

\sus{A transverse slice different from Slodowy's}
Our isomorphisms $\phi$
and $\tphi$ are similar to those conjectured and constructed in \cite{N94, M}.
However, in our case $T_\la$ is \emph{not} the Slodowy's transverse slice but rather
a \emph{different transverse slice} naturally arising from the affine Grassmannian
via the isomorphism $\psi$. In order to illustrate the difference, let us
give an example for $N = 5$ and a nilpotent element $x$ with Jordan blocks of sizes
$3$ and $2$.  If we fix the basis in which the matrix of $x$ has the Jordan canonical form,
i.e.,
\begin{equation}\notag
x = \left (
\begin{matrix}
0 & 1 & 0 & | & 0 & 0\\
0 & 0 & 1 & | & 0 & 0\\
0 & 0 & 0 & | & 0 & 0\\
\hline
0 & 0 & 0 & | & 0 & 1\\
0 & 0 & 0 & | & 0 & 0\\
\end{matrix}
\right ).
\end{equation}
In the Jordan basis the two transverse slices in questions are described by matrices of the form
\begin{equation}\notag
\text{\rm Slodowy's slice} = \left (
\begin{matrix}
a_1 & 1 & 0 & | & 0 & 0\\
a_2 & a_1 & 1 & | & b_1 & 0 \\
a_3 & a_2 & a_1 & | & b_2 & b_1 \\
\hline
c_1 & 0 & 0 & | & d_1 & 1 \\
c_2 & c_1 & 0 & | & d_2 & d_1 \\
\end{matrix}
\right ),
\qquad
\text{\rm our slice} = \left (
\begin{matrix}
0 & 1 & 0 & | & 0 & 0\\
0 & 0 & 1 & | & 0 & 0 \\
a_3 & a_2 & a_1 & | & b_2 & b_1 \\
\hline
0 & 0 & 0 & | & 0 & 1 \\
c_2 & c_1 & 0 & | & d_2 & d_1 \\
\end{matrix}
\right ).
\end{equation}
Let $\{x, h, y\}$ be a Jacobson-Morozov $sl(2)$-triple associated with $x$.
Recall that Slodowy's slice is $x + Z_{\gl(N)}(y)$. Our slice also arises from $\{x, h, y\}$, it can be
described as $x + C\sub \gl(N)$, where $h$ acts on $C$ with non-positive integral eigenvalues,
$C$ is complementary to $[\gl(N), x]$ in $\gl(N)$, and the action of $y$ on $C$ is
\emph{``as close to regular nilpotent as possible"}, cf. \ref{maximalpartition}.
In Slodowy's case the slice is $x + C = x + Z_{\gl(N)}(y)$, so
$h$ acts on $C$ with non-positive integral eigenvalues and
$C$ is complementary to $[\gl(N), x]$ in $\gl(N)$, but by contrast $y$ acts on $C$ by zero.

Our transverse slice is advantageous in the context of this work for three reasons.
First, the isomorphism $\phi$ is given by simple explicit formulas, at least when $c = 0$,
cf. \ref{constructphi} and \cite[3.2]{MVyb},
as opposed to an inductive procedure used in \cite{M}. Second, we are able to
decompose an affine Grassmannian into a disjoint union of quiver varieties, cf. \ref{decompositionAffGrass}.
Finally, our construction provides a natural environment for geometric $(GL(m), GL(n))$ duality,
cf. Section \ref{applicationsRepTheory}.

\sss*{Remark} Slodowy's slice was discovered by Kostant, Peterson and Slodowy,
cf. \cite{Sl, CG} and references therein.

\sus{} The paper is organized as follows. In Section \ref{quiversAsection} we recall some facts on
the quiver varieties of type A. In Section \ref{orbitSection} we recall Grothendieck-Springer-Ginzburg theory
and discuss transverse slices to nilpotent orbits. In Section \ref{grassmannSection} we recall some facts
on Beilinson-Drinfeld Grassmannians and discuss the appearance of our transverse slice in this setting.
Section \ref{mainresults} contains the statement of the Main Theorem and its corollaries. In Section
\ref{quiverConjClassSection} we describe a particular case providing a construction of the conjugacy classes of matrices via quiver varieties. Section \ref{proofMainLemma} contains the proof of the main technical lemma.
Section \ref{proofMainTheorem} finishes the proof of the Main Theorem.
Finally, in Section \ref{applicationsRepTheory} we discuss applications to representation theory.

\sus*{Acknowledgement}
We are grateful to
A. Braverman, I. Frenkel, D. Gaitsgory, V. Ginzburg,
M. Finkelberg, G. Lusztig, A. Maffei, A. Malkin, O. Schiffmann, and W. Wang
for useful discussions,
and to MSRI, IH\' ES and IAS for their
hospitality and support.
The research of I.M. was
supported by NSF.
The research of M.V. was
supported by NSF Postdoctoral
Research Fellowship in 2001-2003.

\se{Quiver varieties of type A}\label{quiversAsection}

\sus{Definitions}

\sss{} Let us consider the Dynkin graph of type $A_{n-1}$ with the
following orientation $\Om$:

\begin{equation}
\nonumber
\xygraph{
[]
!{<0pt,0pt>;<20pt,0pt>:}
*\cir<2pt>{}
!{\save -<0pt,6pt>*\txt{$_1$}  \restore}
-@{{}->} [r]
*\cir<2pt>{}
!{\save -<0pt,6pt>*\txt{$_2$}  \restore}
-@{{}->} [r]
*\cir<2pt>{}
!{\save -<0pt,6pt>*\txt{$_3$}  \restore}
-@{{}->} [r]
*{ \;  \dots \; }
-@{{}->} [r]
*\cir<2pt>{}
!{\save -<0pt,6pt>*\txt{$_{n-2}$}  \restore}
-@{{}->} [r]
*\cir<2pt>{}
!{\save -<0pt,6pt>*\txt{$_{n-1}$}  \restore}
}
\end{equation}

Let $I=\{1,\dots,n-1\}$ be the set of vertices and
$H=\Om\sqcup\overline\Om$
be the set of arrows of our quiver. For an arrow $h\in H$ we denote
by $h'\in I$ its initial vertex and by $h''\in I$ its terminal vertex.

\sss{} Following Nakajima we attach vector spaces $V_i$ and $D_i$
of dimensions $\dim V_i=v_i$ and $\dim D_i=d_i$, $i\in I$ to the
vertices of our quiver i.e. we consider the $I$-graded vector spaces
$V=\oplus_{i\in I}V_i$ and $D=\oplus_{i\in I}D_i$.
Let $v=(v_1,\dots,v_{n-1})$ and
$d=(d_1,\dots,d_{n-1})$ and let $M(v,d)$ be the following affine space:
\[
M(v,w)=\bigoplus_{h\in H}\Hom(V_{h'},V_{h''})\oplus
\bigoplus_{i\in I}\Hom(D_{i},V_{i})\oplus
\bigoplus_{i\in I}\Hom(V_{i},D_{i}).
\]
Following Lusztig \cite{L98} we denote an element in $M(v,w)$ as a triple
$(x,p,q)$, where
\begin{equation}
\begin{split}
x &=(x_h)_{h\in H} \in\bigoplus_{h\in H}\Hom(V_{h'},V_{h''}), \\
p &=(p_i)_{i\in I} \in\bigoplus_{i\in I}\Hom(D_{i},V_{i}), \\
q &=(q_i)_{i\in I} \in\bigoplus_{i\in I}\Hom(V_{i},D_{i}).
\end{split}
\end{equation}

\sss{} In the $A_{n-1}$ case under consideration it is more convenient
to use a different notation.
Following Lusztig and Maffei
we will consider an element in $M(v,w)$ as a quadruple
$(x,\bx,p,q)$. The notation is summarized in the
following diagram:
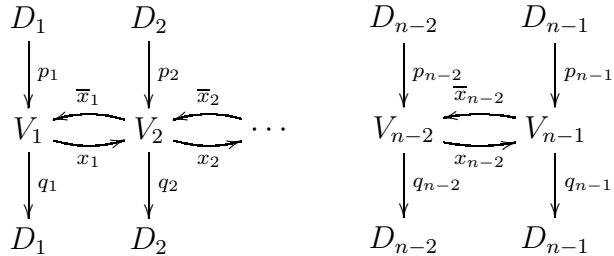
\begin{equation*}
   \xymatrix{
        D_1 \ar[d]^{p_1}  &  D_2 \ar[d]^{p_2}  &  &
        D_{n-2} \ar[d]^{p_{n-2}} & D_{n-1}\ar[d]^{p_{n-1}}\\
        V_1 \ar[d]^{q_1} \ar@/_/[r]_{x_1} &
        V_2 \ar[d]^{q_2} \ar@/_/[r]_{x_2} \ar@/_/[l]_{\bx_1} &
        {\cdots}\ar@/_/[l]_{\bx_2} &
        V_{n-2} \ar[d]^{q_{n-2}} \ar@/_/[r]_{x_{n-2}}  &
        V_{n-1} \ar[d]^{q_{n-1}} \ar@/_/[l]_{\bx_{n-2}}\\
        D_1  &  D_2  &  &  D_{n-2}  &  D_{n-1}
   }
\end{equation*}
Also denote $p_{j \ra i} = \bx_{i} \dots \bx_{j-1} p_j$
and $q_{j \ra i} = q_i x_{i-1} \dots x_{j}$.

The group $G(V)=\prod_{i\in I}GL(V_i)$ acts on $M(v,w)$ in the
following way.
If $g=(g_i)_{i\in I}$ then
\begin{equation}
g(x,\bx,p,q)=(g_{i+1}x_ig_i^{-1}, g_i\bx_ig_{i+1}^{-1},
g_ip_i, q_ig_i^{-1}).
\end{equation}

\sss{} Let us denote by $\mu:M(v,d)\to \mathfrak{g}(V)$ the moment map
associated to this action of $G(V)$. Here $\mathfrak{g}(V)$
is the Lie algebra of $G(V)$. A quadruple $(x,\bx,p,q)$ is in
$\mu^{-1}(c)$, $c=(c_1,\dots,c_{n-1})$
if and only if the following relations are satisfied:
\begin{equation}\label{deformedrelations}
\begin{split}
c_1+\bx_1x_1 &=p_1q_1, \\
c_i+\bx_ix_i &=x_{i-1}\bx_{i-1}+p_iq_i \ \ 2\leq i\leq n-2, \\
c_{n-1} &= x_{n-2}\bx_{n-2}+p_{n-1}q_{n-1}. \\
\end{split}
\end{equation}
We denote the set of all such quadruples by $\La^c(v,d)$.

\sus{A result on invariant polynomials} Following \cite{L98}
let $\RR$ be the algebra of regular functions $M(v,d)\to \C$
and let $\RR(\La)$ be the algebra of regular functions
$\La^c(v,d)\to \C$. The action of $G(V)$ on $M(v,d)$ (resp. $\La^c(v,d)$)
induces an action of $G(V)$ on $\RR$ (resp. $\RR(\La)$).
Following Lusztig \cite[1.2]{L98} we describe two groups
of invariant polynomials in $\RR^{G(V)}$.

(a) Let $h_1,h_2\dots,h_r$ be a cycle in our graph, that is a sequence
in $H$ such that $h_1''=h_2', h_2''=h_3',\dots,h_{r}''=h_1'$.
This cycle defines a $G(V)$-invariant polynomial in $\RR^{G(V)}$
given by
$(x,p,q)\mapsto \Tr(x_{h_r}x_{h_{r-1}}\dots x_{h_1}):
V_{h_1'}\to V_{h_1'}$.

(b) Let $h_1,h_2\dots,h_r$ be a path in our graph, that is a sequence
in $H$ such that $h_1''=h_2', h_2''=h_3',\dots,h_{r-1}''=h_r'$.
This path together with a linear form $\chi$ on
$\Hom(D_{h_1'}, D_{h_r''})$ defines a $G(V)$-invariant polynomial in
$\RR^{G(V)}$ given by
$(x,p,q)\mapsto \chi(q_{h_r''}x_{h_r}x_{h_{r-1}}\dots x_{h_1}p_{h_1'})$.

\theo\label{LusztigTheorem} \cite[Theorem 1.3, 5.8]{L98} The algebra $\RR(\La)^{G(V)}$
is generated by the invariant polynomials of types (a) and (b) above
for $(x,p,q)\in\La^c(v,d)$.

Following \cite{M}, in the $A_{n-1}$ case we can improve the above
theorem as follows. We switch back to Maffei's notation.

\slemm{} Let  $h_1,h_2\dots,h_r$ be a cycle in our quiver. Then
$$
\Tr(x_{h_r}x_{h_{r-1}}\dots x_{h_1})=
\Tr(P),
$$
where $P$ is some polynomial of $q_{l\ra j}p_{i\ra l}$,
$i,j\in\{1,\dots,n-1\}$ (necessarily $l\leq\min(i,j)$).

\begin{proof} Easily follows from relations (\ref{deformedrelations}).
\end{proof}

\slemm{} Let $h_1,h_2\dots,h_r$ be a path in our graph
and let $\chi$ be a linear form on\linebreak
$\Hom(D_{h_1'}, D_{h_r''})$. Then
$$
\chi(q_{h_r''}x_{h_r}x_{h_{r-1}}\dots x_{h_1}p_{h_1'}),
=\chi(P),
$$
where $P$ is some polynomial of $q_{l\ra j}p_{i\ra l}$,
$i,j\in\{1,\dots,n-1\}$ (necessarily $l\leq\min(i,j)$).

\begin{proof} Easily follows from relations (\ref{deformedrelations}).
\end{proof}

Notice that $\Tr:D_i\to D_i$ is a linear form on $\Hom(D_i,D_i)$.
Now the Lusztig's theorem \ref{LusztigTheorem} and the lemmas above imply
the following.

\theo{}\label{invtheorem} The algebra of invariant functions
$\RR(\La)^{G(V)}$
is generated by the invariant polynomials
$\chi(q_{l\ra j}p_{i\ra l})$, where
$i,j\in\{1,\dots,n-1\}$, $1\leq l\leq\min(i,j)$,
and $\chi$ is a linear form on
$\Hom(D_i, D_j)$.

\sss{} Following Nakajima \cite{N98} and Lusztig \cite[2.11]{L98}
we say that a quadruple
$(x,\bx,p,q)$ is stable if for any $I$-graded subspace $U$ of $V$
containing
$\Im p$ and preserved by
$x$ and $\bx$, we have $U=V$. The set of all stable quadruples in
$\La^c(v,d)$ is denoted by
$\La^c_s(v,d)$.

The following easy lemma is lifted from Maffei, \cite[Lemma 14]{M}.

\slemm{}\label{stablelemma}
If $(x,\bx,p,q)\in\La^c(v,d)$ then $(x,\bx,p,q)$ is stable if
and only if
for all $1\leq i\leq n-1$
\begin{equation}
\Im\ x_{i-1}+\sum_{j=i}^{n-1} \Im\ p_{j\ra i}=V_i.
\end{equation}

\sus{Nakajima's quiver variety \cite[3.12]{N98}}\label{quiverdefs}
The quiver variety $\fM(v,d)$
is the geometric quotient of $\La^c_s(v,d)$ by $G(V)$.
In particular the set of geometric points of $\fM$
is  $\La^c_s(v,d)/G(V)$.
Below we only consider such $(v,d)$ that $\fM(v,d)$ is nonempty,
see \cite[10]{N98}, \cite[Lemma 7]{M} for explicit conditions on $(v,d)$.

We can also consider the affine algebro-geometric quotient of
$\La^c(v,d)$ by $G(V)$, which we denote by
\begin{equation}
\fM_0=\La^c(v,d)//G(V)=\Spec \RR(\La^c(v,d))^{G(V)}.
\end{equation}

We have a natural map $p:\fM(v,d)\ra \fM_0(v,d)$.
Following Maffei we denote
$$
\Im\ p=\fM_1(v,d)\subset\fM_0(v,d).
$$
Finally, let  $\fL(v,d):= p^{-1}(0)\sub \fM(v,d)$ and
denote by
$
\HH(\fL(v,d))
$
its  top-dimensional Borel-Moore homology.

\sus{$SL(n)$-modules}\label{slnmodules} In this subsection $c=0$.

\stheo{}
\label{theonakajima}
\cite[\emph{10.ii}]{N98}
The space $\pl_v \HH(\fL(v,d))$ has the structure of a
simple $SL(n)$-module $W_d$ with the highest weight $d$
(i.e., $\sum_I\ d_i\om_i$ for the  fundamental weights
$\om_i$).
The summand
$\HH(\fL(v,d))$ is the weight space for
the weight
$d-Cv$, where $C$ is the Cartan matrix of type $A_{n-1}$.

In particular, the module $W_d$ has a basis arising from
the irreducible components of $p^{-1}(0)$,
or more precisely the weight space
$W_d(d-Cv)$ has a basis indexed by $\Irr\fL(v,d)$.
Following Lusztig \cite{L00}, we call this basis \emph{semicanonical}.

\sss{From $SL(n)$ to $GL(n)$}\label{fromSLntoGLn}
\label{sltogl}
We may consider $\pl_v \HH(\fL(v,d))$ as a representation $W_{\cla}$
of $GL(n)$
with highest weight $\cla$,
where $\cla=\cla(d)=(\cla_1, \cla_2,\dots,\cla_n)$
is a partition of
$N=\sum_{j=1}^{n-1}jd_j$ defined as follows:
$\cla_i=\sum_{j=i}^{n}d_j$ (here $d_n=0$). Then
$\HH(\fL(v,d))$ is the weight space $W_{\cla}(a)$,
where $a_i=v_{n-1}+\sum_{j=i}^{n}(d-Cv)_j$ (here $(d-Cv)_n=0$),
cf. \cite[8.3]{N94}.

\se{Grothendieck-Springer-Ginzburg
theory and conjugacy classes of matrices}\label{orbitSection}

In this section we fix a vector space $D$ of dimension $N$.

\sus{Definitions of bases}\label{definitionBases}

\sss{} Let $\NN=\NN(D)$ be the nilpotent cone in $\End(D)$.
Let $a=(a_1,\dots,a_n)$ be a $n$-tuple of integers such that
$N=\sum_{i=1}^{n}a_i$.
We denote the variety if $n$-step flags in $D$ and its
connected components as follows:
\begin{equation}
\begin{aligned}
\FF^n &=\{0=F_0\subseteq F_1\subseteq F_2\subseteq \dots\subseteq F_n=D\}, \\
\FF^{n,a} &=\{0=F_0\subseteq F_1\subseteq F_2\subseteq
\dots\subseteq F_n=D ~|~\dim F_i-\dim F_{i-1}=a_i\}.
\end{aligned}
\end{equation}
It is well known that we have the following description of the
cotangent bundle $\TN^n=T^*\FF^n$ to this flag variety
and its connected components
\begin{equation}
\TN^{n,a}=T^*\FF^{n,a}=\{(x,F)\in\NN\times\FF^{n,a} ~|~ x(F_i)
\subseteq F_{i-1} \}.
\end{equation}

Denote by $\mmm:\TN^{n}\ra \NN$ the projection onto the first factor,
and by $\mmm_a$ the restriction of $\mmm$ to $\TN^{n,a}$.

\sss{}\label{ginzburgfiber}
Let $\cla=\cla_1\geq\dots\geq\cla_n$,
$N=\sum_{i=1}^{n}\cla_i$ be a partition of $N$ and let
$\la=(\la_1,\dots,\la_m)$, be the dual partition.
Let $x\in\NN$ be a nilpotent element of type $\la$,
that is, $x$ has Jordan blocks of sizes
$\la_1, \dots, \la_m$.
We will denote the fiber $\mmm^{-1}(x)$ by $\FF^n_x$
and its connected components $\mmm^{-1}(x)\cap \FF^{n,a}$
by $\FF^{n,a}_x$.

\sss{} Let us extend the picture above as follows.
Consider the following subbundle of the trivial
vector bundle $\\gl(D)\times\FF^n$ (resp. $\gl(D)\times\FF^{n,a}$):
\begin{equation}
\begin{aligned}
\tfg=\tfg^n= & \{(x,F)\in\gl(D)\times\FF^n ~|~ x(F_i)
\sub F_i \}, \\
\tfg^{n, a}= & \{(x,F)\in\gl(D)\times\FF^{n,a} ~|~ x(F_i)
\sub F_i \}.
\end{aligned}
\end{equation}
We will denote the projection to the first factor
by $\tmm: \tfg\to \fg=\gl(D)$.
More notation: $\tfg_x:=\tmm^{-1}(x)$
and $\tfg^{n,a}_x:=\tmm^{-1}(x)\cap \tfg^{n, a}$.

\sss{}\label{conjclass}
Let us fix $x\in\End(D)$ with the spectrum
(=set of eigenvalues)
$E\sub\ \A^1$
such that $|E|\leq n$. For $e\in E$ let
the restriction of $(x-e\Id_D)$
to the generalized $e$-eigenspace of $x$ be a nilpotent
of type $\mu(e)$ where
$\mu(e)=(\mu_1(e)\geq \mu_2(e)\geq \dots \geq \mu_{m(e)}(e))$
is a partition and $|\mu(e)|=\sum_{i=1}^{m(e)} \mu_i(e)=l(e)$,
so $l(e)$ is the multiplicity of $e$.
For every partition $\mu(e)$ consider its dual
$\cmu(e)=(\cmu_1(e)\geq \cmu_2(e)\geq \dots \geq \cmu_{n(e)}(e))$,
so $n(e)$ is the size of the largest Jordan block associated with $e$.
Let $\tmu=\{\mu(e)\}_{e\in E}$ be the collection of partitions
for all eigenvalues of $x$.

The data $E, \tmu$ define the conjugacy class of $x$ (Jordan canonical form).
Let us denote this conjugacy class by $\OO_{E,\tmu}$.

Let us assume now that $\sum_{e\in E} n(e)=n$.
Then the set of pairs
\begin{equation}\label{setM}
M=\{(e,\cmu_i(e))~|~ e\in E, 1\leq i\leq n(e)\} \subset E \times \Z^n
\end{equation}
is an $n$-element subset of $E \times \Z^n$. Let us take an arbitrary bijection
$\be:[1,n]\to M$, where $[1,n]$ is the set of integers from $1$ to $n$.
Let $\be_1:[1,n]\overset{\be}\ra M\to E$
be the composition of $\be$ with the projection of $M$ to the first
factor, and let $\be_2:[1,n]\overset{\be}\ra M\hookrightarrow E \times \Z^n\to \Z^n$
be the composition of $\be$
with the inclusion of $M$ into $E \times \Z^n$ and the projection to the second factor.
Denote $a=(a_1,\dots, a_n) = (\be_2(1), \dots, \be_2(n))$.

Now we can consider
\begin{equation}\label{gnabth}
\tfg^{n,a,E,\tmu}=\{(x,F)\in\barr\OO_{E,\tmu}\times\FF^{n,a} ~|~ x(F_i)
\sub F_i
\text{ and } x \text{ acts on } F_l/F_{l-1}
\text{ as } \be_1(l)\Id \}.
\end{equation}
We will still denote the projection to the first factor
by $\tmm: \tfg^{n,a,E,\tmu}\to \barr\OO_{E,\tmu}$.
Now we need the following.

\slemm{} The variety $\tfg^{n,a,E,\tmu}$ is smooth and  connected,
the map $\tmm$ is projective and
$$
\dim\tfg^{n,a,E,\tmu}=\dim \OO_{E,\tmu}=
N^2-\sum_{e\in E}\sum_{i\in [1,n(e)]} \cmu^2_i(e).
$$

\begin{proof} Actually $\tfg^{n,a,E,\tmu}$ is a vector bundle
over $\FF^{n,a}$ with the fiber over a particular flag $F$
being $P(F)/L(F)$ where $P(F)$ is the parabolic preserving $F$
and $L(F)$ its Levi factor. Also, if $x\in\OO_{E,\tmu}$, then
$\tmm^{-1}(x)
\cap \tfg^{n,a,E,\tmu}$ is a point.
In fact, the conjugacy class $\OO_{E,\tmu}$ is a deformation
of the nilpotent class $\OO_\mu$ where $\cmu$ is the partition
obtained from the $n$-tuple $a=(a_1,\dots, a_n)$ as above by ordering the elements
in the non-increasing order. The variety $\tfg^{n,a,E,\tmu}$
is isomorphic to $\TN^{n,a}$. In particular, $\dim \OO_{E,\tmu}=\dim \OO_\mu$.
\end{proof}

\sss{} For a finite dimensional algebraic variety $X$ we denote
by $H(X)$ its top-dimensional Borel-Moore homology $H^{\rm BM}_{\dim X}(X)$.
In particular, we denote
$$
\begin{aligned}
H(\FF^n_x)& :=\bigoplus_{a}H^{\rm BM}_{\dim\FF^{n,a}_x}(\FF^{n,a}_x), \\
H(\tfg_x)& :=\bigoplus_{a}H^{\rm BM}_{\dim\tfg^{n,a}_x}(\tfg^{n,a}_x).
\end{aligned}
$$

The following theorem is due to Ginzburg and Braverman-Gaitsgory.

\theo{}\label{theoginzburg}
\begin{enumerate}
\item \cite[4.2]{CG}
Let $\check x$ be a nilpotent of type $\cla$.
The space $H(\FF^n_{\check x})$ has the structure of a
$\gl(n)$-module $W_{\la}$ with the highest weight $\la$.
\item \cite{BG}
Let $x$ be a nilpotent of type $\la$.
The space $H(\tfg_x)$ has the structure of
a $\gl(n)$-module $W_{\la}$ with the highest weight $\la$.
\end{enumerate}

In particular, the module $W_{\la}$ has two bases:
\begin{enumerate}
\item A basis indexed by
$\Irr\FF^n_{\check x}$. More precisely, the weight space $W_{\la}(a)$
has a basis indexed by $\Irr\FF^{n,a}_{\check x}$.
It was shown in \cite{Sav} that this basis
coincides with the semicanonical basis defined in \ref{slnmodules}.
\item A basis indexed by $\Irr\tfg_x$ (relevant irreducible components).
More precisely, the weight space $W_{\la}(a)$
has a basis indexed
by $\Irr\tfg^{n,a}_x$. We call this the \emph{Spaltenstein basis}.
\end{enumerate}

\sss{Remark} It was established in \cite{BGV} that
the Spaltenstein basis as above coincides with the
Mirkovi\' c-Vilonen basis of \cite{MV}. As far as we know the question about
the relationship between the semicanonical (as well as Lusztig's canonical \cite{L90, L91})
and Mirkovi\' c-Vilonen bases remains open.

\sus{On normal (transverse) slices}\label{ontransverseslices}

Let $\fg=\gl(D)$ and $G = GL(D)$.

\sss{
Normal slices to nilpotent orbits
}

We will say that a normal slice (in $\fg$) to
a nilpotent orbit $\al$ at $e\in\al$,
is a submanifold $S$ of $\fg$ such that
\begin{enumerate}
\item
(Infinitesimal normality.)
$T_e\al\ \pl\  T_eS=\ \fg$,
(cf. \cite[3.2.19]{CG}) and
\item (Contraction.)
There is an action of $G_m$ on $S$ which
contracts it to $e$ and preserves intersections with
the Lusztig strata
in $\fg$. (For the definition of Lusztig strata cf. \cite[5.5]{Mir} and references therein.)
\end{enumerate}
We will use the terminology "normal slice" and "transverse slice" interchangeably.

\lemm
For a normal slice $S$
\begin{enumerate}
\item $S\cap \al=\{e\}$.
\item $S$ meets Lusztig stratum $\be$ iff $\al\sub \barr\be$.
\item $S$ meets Lusztig strata transversally.
\end{enumerate}

\lemm\label{sufficientnormaldata}
A sufficient data for a normal slice to the orbit $\ ^Ge$ at $e$
is given by a pair $(h,C)$
where $h\in\fg$
is semisimple
integral (i.e., eigenvalues of $\ad\ h$ are integral),
and $[h,e]=2e$;
while
$C\sub \fg$ is an $h$-invariant vector subspace
complementary to $T_e(\al)=[\fg,e]$,
such that the eigenvalues of $h$ in $C$ are
${\le 1}$. Then $S=e+C$ is a normal slice.

\begin{proof}
Such $h$ lifts to a homomorphism $\io: G_m\to G$ and we can construct
an action of $G_m$ on the vector space $\fg$ by
$
s\ast x=\ s^{-2}\ \cd\ ^{\io(s)}x, \ s\in G_m,\ x\in\fg
$
which fixes $e$ and preserves $e+C$.
\end{proof}

\sss{} For a nilpotent $e$ let $\{e, h, f\}$ be a Jacobson-Morozov $sl(2)$-triple.
We can build normal slices to the nilpotent orbit $^Ge$ at $e$ using $h$ and $f$.

\sss{Example: Slodowy's slice} Take $h, f$ from a Jacobson-Morozov $sl(2)$-triple,
and let $C = Z_{\fg} (f)$. Clearly, the conditions of Lemma \ref{sufficientnormaldata}
are satisfied, and $S = e + Z_{\fg} (f)$ is the best-known example of a normal slice.

\sss{Another slice}\label{basisfreeslice}
We will consider another slice arising from
a Jacobson-Morozov $sl(2)$-triple $\{e, h, f\}$.
First, let ${}^h\fg_{\leq 0} \sub \fg$ be the $\{h, f\}$-invariant subspace
such that the eigenvalues of $h$ in ${}^h\fg_{\leq 0}$ are ${\le 0}$. Then
$f$ will act as a nilpotent in ${}^h\fg_{\leq 0}$ and to build a normal slice $S = e + C$
if suffices to choose $C\sub {}^h\fg_{\leq 0}$ complementary to $T_e(\al)=[\fg,e]$.

If we choose $C = \ker_{\fg} (f)\sub {}^h\fg_{\leq 0}$ we recover the Slodowy's slice.

We would like to consider a $C \sub {}^h\fg_{\leq 0}$ with the property that
$f$ restricted to $C$ is \emph{``as close to regular nilpotent as possible"},
cf. \ref{maximalpartition} for more details.
In particular, if $e$ is regular, then $f$ restricted to our $C$ will be regular.

More precisely, the vector space $D$ considered
as an $sl(2)$-module decomposes as:
\begin{equation}
D = \bigoplus_i M_i\otimes L_i, \qquad{ \rm{\text{ and }}}\qquad \End(D) = D^{\ast} \otimes D \simeq \bigoplus_{i, j} \Hom(M_j, M_i)\otimes L_j^{\ast} \otimes L_i,
\end{equation}
where $L_i$ is a simple $sl(2)$-module of highest weight $i$, $\dim L_i = i+1$, and $M_i$ is its multiplicity in the decomposition above.

Now consider $C$ to be a subspace
\begin{equation}
C = \bigoplus_{i, j} \Hom(M_j, M_i)\otimes \ker_{L_j^{\ast}}(f^{i+1})
\otimes \ker_{L_i}(f) \subseteq \End(D),
\end{equation}
where $\ker_{L_j^{\ast}}(f^{i+1})$ (resp. $\ker_{L_i}(f)$) is the kernel of the natural action of
$f^{i+1}$ (resp. $f$) on $L_j^{\ast}$ (resp. $L_i$.). Notice that $\dim \ker_{L_j^{\ast}}(f^{i+1}) = i+1$,
and $\dim \ker_{L_i}(f) = 1$, and also $\dim C = \dim Z_{\fg}(f)$.

It is elementary to see that $h, C$ satisfy the conditions of Lemma \ref{sufficientnormaldata}, and thus
$S = e + C$ is a normal slice.

\sss{``As close to regular nilpotent as possible"}\label{maximalpartition}
Let $e$ be acting on $D$ as a nilpotent of type $\la=(\la_1\geq\la_2\geq\dots\geq \la_m)$.
Then we could say that
\begin{equation}
D = \bigoplus_{i=1}^m  L_{\la_i-1}, \qquad{ \rm{\text{ and }}}\qquad \End(D) = \bigoplus_{i, j = 1}^m
L_{\la_j-1}^{\ast} \otimes L_{\la_i-1},
\end{equation}
where $L_{\la_i-1}$ is a simple $sl(2)$-module of highest weight $\la_i - 1$, $\dim L_i = \la_i$.
If $\la_i \geq \la_j$ we have
$$
L_{\la_j-1}^{\ast} \otimes L_{\la_i-1} = L_{\la_i + \la_j - 2} \oplus L_{\la_i + \la_j - 4}
\oplus \dots \oplus L_{\la_i - \la_j},
$$
$\la_j$ summands in all. Let $f$ be the element of a
Jacobson-Morozov $sl(2)$-triple. Observe that $f$ restricted to $C\cap (L_{\la_j-1}^{\ast} \otimes L_{\la_i-1})$
acts as a regular nilpotent. It is easy to see
that $f$ acts on $C$ defined as above as a nilpotent of type
\begin{equation}
\la_f= (\la_1 \geq \la_2\geq \la_2\geq \la_2 \geq \dots \geq \la_m \geq \dots \geq \la_m),
\end{equation}
where $\la_k$ repeats with multiplicity $2k - 1$, $1\leq k \leq m$. Then $\la_f$ is a partition of $\sum_{i=1}^{m}(\cla_i)^2 = N^2 - \dim\OO_{\la}$ and the largest such partition
possible for $C \sub  {}^h\fg_{\leq 0}$ and $C$ being complementary to $T_e(\al)=[\fg,e]$.
By contrast in the Slodowy's situation $f$ acts
on $C = \ker_{\fg}(f)$ as $0$ and so its type $(1, \dots, 1)$ is the smallest possible partition of $\sum_{i=1}^{m}(\cla_i)^2$.

\sus{Our slice in Jordan basis}\label{sliceJordanBasis}

We will adjust the notation a bit here: the nilpotent $e$ will be denoted $x$ in this subsection.

\sss{}\label{nilpnormalslice}
Again, let $D$ be a vector space, $\dim D=N$,
and $\NN$ be the nilpotent cone in $\End(D)$.
Let $x$ be
a nilpotent operator of type $\la=(\la_1\geq\la_2\geq\dots\geq \la_m)$.
Moreover, $e_{k,i}$, $1\leq k\leq \la_i$
be a basis in $D$ in which $x$ is exactly the direct
sum of nilpotent blocks and $x$ restricted to the span of
$\{e_{k,i}\ | \ 1\leq k\leq \la_i \}$ is the Jordan block of size $\la_i$,
that is $x :e_{k,i}\mapsto e_{k-1,i},\ e_{1,i}\mapsto 0$.

Now define:
\begin{equation}
T_x:=\{ x+f,\  f\in\End(D) \ |\
f^{l,j}_{k,i} =0, \text { if } k\neq \la_i,
\text{ and }
f^{l,j}_{\la_i,i} =0, \text { if } l>\la_i \},
\end{equation}
where $f^{l,j}_{k,i}:\C e_{l,j}\to\C e_{k,i}$
are the matrix elements of $f$ in our basis. For example,
if $\la=(\la_1\geq\la_2) = (3, 2)$ the matrices in $T_x$ in the
basis $e_{k,i}$, $1\leq k\leq \la_i$ will have the form
\begin{equation}\label{tsliceexample}
\left(
\begin{matrix}
0 & 1 & 0 & | & 0 & 0 \\
0 & 0 & 1 & | & 0 & 0 \\
f^{1, 1}_{3, 1} & f^{2, 1}_{3, 1} & f^{3, 1}_{3, 1} & |& f^{1, 2}_{3, 1} & f^{2, 2}_{3, 1}\\
\hline
0 & 0 & 0 & | & 0 & 1 \\
f^{1, 1}_{2, 2} & f^{2, 1}_{2, 2} & 0 & | & f^{1, 2}_{2, 2} & f^{2, 2}_{1, 2} \\
\end{matrix}
\right ).
\end{equation}

The set $T_x$ (denoted by $e + C$ above) will sometimes be denoted by $T_\la$.

\sss{}
For $\mu$ such that $\OO_{\la} \sub \barr\OO_{\mu}$ define
$$
T_{x,\mu}:=T_x\cap \barr\OO_{\mu}.
$$
We have seen in \ref{basisfreeslice} that

\slemm{}
$T_x$ is a transverse slice to the orbit of $x$.
In particular,
\begin{equation}
\dim T_{x,\mu}=\dim \OO_{\mu}-\dim\OO_{\la}=
\sum_{i=1}^{\la_1}(\cla_i)^2-\sum_{i=1}^{\mu_1}(\cmu_i)^2.
\end{equation}

\sss{}\label{tildeta}
Let $\la\leq \mu$ and take
any permutation $a=(a_1,\dots,a_n)$ of the dual partition $\cmu$.
We will restrict the resolution
$\mmm$ to the slice $T_{x,\mu}$ :
$$
\TT_{x}^{a}:=\mmm_{a}^{-1}(T_{x,\mu})\subset
\TN^{n,a}.
$$

\slemm{} The variety $\TT_{x}^{a}$ is smooth and connected
of dimension $\sum_{i=1}^{\la_1}(\cla_i)^2-\sum_{i=1}^{\mu_1}(\cmu_i)^2$.
It is nonempty if and only if $x\in\barr\OO_{\mu}$.

The map $\mmm_a: \TT_{x}^{a}\to T_{x}\cap \barr\OO_{\mu}$
is projective.

\begin{proof}
$\TT_{x}^{a}$ is smooth because $G\cd T_{x,\mu}$ is open in $\fg$
and near $T_{x,\mu}$ it is a product of
$T_{x,\mu}$ and the orbit $G\cd x$.
The dimension counts follow from
\begin{equation}\nonumber
\dim\OO_{\la}=N^2-\sum_{i=1}^{\la_1}(\cla_i)^2=N^2-\sum_{i=1}^{m}(2i-1)\la_i.
\end{equation}
\end{proof}

\sss{} We also need to study the intersection $T_x\cap\OO_{E,\tmu}$,
where $\OO_{E,\tmu}\sub \End(D)$ is a conjugacy class defined in
\ref{conjclass}.

\slemm{} We have
$$
\dim T_x\cap\OO_{E,\tmu}=
\sum_{i=1}^{\la_1}(\cla_i)^2-\sum_{e\in E}\sum_{i\in [1,n(e)]} \cmu^2_i(e).
$$
Moreover, $T_x\cap\OO_{E,\tmu}$ is nonempty if and only of
$x\in\barr\OO_{\mu}$, where $\mu$ is obtained from $\tmu$ as in
\ref{conjclass}.

\begin{proof} Follows from general smoothness results.
\end{proof}

\slemm{} The variety
$\tmm(T_x\cap\OO_{E,\tmu})\cap \tfg^{n,a,E,\tmu}$ is smooth
and connected of dimension equal to $\dim T_x\cap\OO_{E,\tmu}$.

\begin{proof}
The proof is the same as above, for connectedness cf. \cite{Sp}.
\end{proof}

\se{Beilinson-Drinfeld Grassmannians of type A}\label{grassmannSection}

We recall some standard facts about the affine Grassmannians
of type A. In this section $G=GL(m)$ unless indicated otherwise.

\sus{Local picture}

\sss{} Let $m$ be a positive natural number, and $V$ a vector
space of dimension $m$. Let us fix a direct sum decomposition of $V$
\begin{equation}
V=V_1\oplus\dots\oplus V_m,
\end{equation}
where $\dim V_i=1$, $1\leq i\leq m$. Let us fix nonzero elements
$\bfe_i\in V_i$. The set $\{\bfe_1,\dots,\bfe_m\}$ is a basis in $V$.

Let $O:=\C[[z]]$ be the ring of formal power series in $z$
and $K:=\C((z))$ be its field of fractions. Let $V(K)=V\otimes K$
and let $L_0=V\otimes O$.
A \emph{lattice} $L$ in $V((z))$
is an $O$-submodule of $V(K)$
such that $L\otimes_{O} K=V(K)$.

The affine Grassmannian $\GG_{G}$ is a (reduced)  ind-scheme
whose $\C$-points can be described as all lattices in $V(K)$
or as $G(K)/G(O)$.
Its connected components $\GG_{(N)}$
are indexed by integers $N\in \Z$. If $N\ge 0$
then $\GG_{(N)}$ contains the finite dimensional subscheme
\begin{equation}\label{ggn}
\GG_N=\{\text{lattices } L \text{ in } V((z)) \text{ such that }
L_0\subseteq L, \dim L/L_0=N \}.
\end{equation}
To a dominant coweight $\la\in\Z^m$ of $G$, one attaches
the lattice $L_{\la}=\ \pl_1^m\ \C[[z]]\cd z^{-\la_i}e_i$.

The $G(O)$-orbits $\GG_\la$ in $ \GG_{G}$
are parameterized by the dominant coweights (partitions) $\la$
via $\GG_\la=G(O)\cd L_\la$.

The $G(O)$-orbits in $\GG_{N}$ correspond to
partitions
$\mu=(\mu_1\geq\mu_2\geq\dots\geq\mu_m)$
of $N$ into at most $m$ parts.
These orbits can be explicitly
described as follows:
\begin{equation}
\GG_{\mu}=\{L\in\GG_N ~|~z \text{ restricted to }
L/L_0 \text{ has Jordan blocks of sizes } \mu_1,\mu_2,\dots,\mu_m \}.
\end{equation}

\sss{} Let $G=PGL(m)$. Then the points of $\GG_G$ can be thought of as
lattices in $V((z))$ only up to a shift by $z$, or as
$PGL(m,K)/PGL(m,O)$. Set theoretically
$\GG_{PGL(m)}$ is a union of $m$ connected components of $\GG_{GL(m)}$.

%\sss{Remark} In fact since we are talking about $\GG_{PGL(m)}$
%we need only lattices up to a shift by $z$, so we fix a representative
%$L$ in each equivalence class such that $z^{-1}L_0\sub L$
%but $z^{-2}L_0\not\sub L$.

\sss{} The orbits of $PGL(m,O)$ on $\GG_{PGL(m)}$ are parametrized
by the dominant weights of the Langlands dual group
$^{L}PGL(m)=SL(m)$.
If we consider $\mu=(\mu_1\geq\mu_2\geq\dots\geq\mu_m)$
defined up to simultaneous shift of by an integer
as a dominant weight of $SL(m)$ then
the $PGL(m,O)$-orbit $\GG_{\mu}$ is described as follows:
\begin{equation}
\GG_{\mu}=\{L\in\GG_N ~|~z \text{ restricted to }
L/L_0 \text{ has Jordan blocks of sizes } \mu \}.
\end{equation}
This is well defined since the lattice $L$ is considered up to a shift
by $z$.

\sus{Global picture}\label{globalpicture}

Let $X$ be a curve, which in our case will always be $\A^1$.
Let $\A^{(n)}=\A^1\times \dots \times \A^1//{\mathfrak S}_n$
be the symmetric $n$-fold product of $\A^1$

Beilinson-Drinfeld Grassmannian \cite{BD, MV, MV2} is
a (reduced) ind-scheme $\BDG_{\A^{(n)}}$ whose $\C$-points
are described as follows:
\begin{equation}
\BDG_{\A^{(n)}}(\C)=\{(b,\VV, t)~ |~ t: \VV_{X-E}\to (X\times V)|_{X-E}
\text{ is an isomorphism}~\},
\end{equation}
where $b=(b_1,\dots, b_n)\in \A^{(n)}$, $E=\{b_1,\dots, b_n\}\sub \A^1$,
$\VV$ is a vector bundle of rank $m$, and $t$ is the trivialization
of $\VV$ off $E$. The pairs $(\VV,t)$ are considered up to an
isomorphism. If we fix $b=(b_1,\dots, b_n)$
(and therefore $E=\{b_1,\dots, b_n\}$)
then the corresponding
ind-subscheme of $\BDG_{\A^{(n)}}$ is called the fiber
of $\BDG_{\A^{(n)}}$ at $b$ and is denoted by $\GG^{BD}_b$.
If $n=1$ we will also write $\GG_{e}$ for $e\in \A^1$.
It is well known \cite{BD, MV} that
\begin{equation}\label{bdproduct}
\BDG_b=\prod_{e\in E} \GG_e .
\end{equation}

\sss{}\label{bdgpolynomial}
Let $\C[z]$ be the ring of polynomials in $z$
and $\C(z)$ be its field of fractions i.e. rational functions.
Let $V(z)=V\otimes \C(z)$
and let $\LL_0=V\otimes O$.
A \emph{lattice} in $V(z)$
is an $\C[z]$-submodule $\LL$ of $V(z)$
such that $L\otimes_{\C[z]} \C(z)=V(z)$.

The points of $\GG_b$ can be described as lattices $\LL$ in
$V(z)=V\otimes \C(z)$
such that their localizations $L(e)$ at $e\in \A^1-E$
are isomorphic to $L_0(e)=V\otimes \C[[z-e]]$.
Define:
$$
\BDG_N=\{ \text{ lattices } \LL\supseteq  \LL_0 ~
|~ \dim \LL/\LL_0=N~ \}.
$$

Slightly generalizing the exposition \cite[Partie I]{Ngo},
we fix a polynomial $P$ of degree $n$, where
$n\leq N\leq mn$.
Define
$$
\BDG_N(P)=\{ \text{ lattices } \LL\supseteq  \LL_0 ~
|~ \dim \LL/\LL_0=N  \text{ and } P(z|_{\LL/\LL_0})=0 \} ,
$$
where $z|_{\LL/\LL_0}$ is the linear operator
on $\LL/\LL_0$ obtained by the restriction of $z$.

Let $P=\prod_{e\in E}(z-e)^{n(e)}$. Then a version of
(\ref{bdproduct}) is
\begin{equation}\label{polybdg}
\BDG_N(P)=\bigsqcup_{
\substack {l(e)\geq n(e) \\
\sum_{e\in E} l(e)=N}}
\prod_{e\in E} (\GG_e)_{l(e)},
\end{equation}
where the finite dimensional subscheme
$(\GG_e)_{l(e)}$ of the affine Grassmannian $\GG_e$
is defined as in (\ref{ggn}).

Finally, if $(b_1, \dots, b_n)\in \A^{(n)}$ and
$E=\{b_1, \dots, b_n\}\subset \A^1$, then
we set $\BDG_{N,b}(P):=\BDG_b\cap \BDG_N(P)$.

\sss{}\label{convolutionglobal}
Let $a=(a_1,\dots, a_n)$ such that $\sum_{i=1}^n a_i=N$.
Let us introduce a convolution Grassmannian
$\TBDG^{n,a}_N$ as the (reduced) scheme whose $\C$-points are
$n$-step flags of lattices in $V(z)$:
$$
\TBDG^{n,a}_N
=\{ \LL_0\sub \LL_1\sub \dots \sub \LL_n=\LL ~|~
\dim\LL_i/\LL_{i-1}=a_i
\text{ for }
1\leq i\leq n  \} ,
$$
where $\LL_0=V\otimes \C[z]$.
We have a map $\pi^{n,a}_N=\pi: \TBDG_N^{n,a}\to \BDG_{N}$ such that
$\pi: (\LL_0\sub \LL_1\sub \dots \sub \LL_n)\mapsto \LL=\LL_n$.

Let $(b_1, \dots, b_n)\in \A^{n}$.
Let us also introduce a subscheme in the fiber of $\TBDG^{n,a}$
over the point $(b_1, \dots, b_n)\in \A^{(n)}$.

$$
\TBDG^{n,a}_b
=\{ (\LL_0\sub \LL_1\sub \dots \sub \LL_n)\in \TBDG^{n,a}_N~|~
\LL_n\in\BDG_b,
\text{ and } z \text { acts on }
\LL_i/\LL_{i-1} \text{ as } b_i \} ,
$$

Finally, if $\{b_1, \dots, b_n\}=E\subset \A^1$, and
$P$ is a polynomial as in (\ref{polybdg}), then we define
$\TBDG^{n,a}_b(P)=\TBDG^{n,a}_b\cap \pi^{-1}(\BDG_N(P))$.

\sss{}\label{convolutionlocal}
Let us also consider the local version of the convolution
Grassmannian. Let $\mu$ be a partition of $N$ into at most $m$ parts
and let $\GG_\mu\sub \GG_N$ be a $G(O)$-orbit in $\GG_N$.
Let $a=(a_1,\dots, a_n)$ be a permutation of the dual partition $\cmu$.
Consider the
(reduced) ind-scheme
$\TGG^a_{\mu}=\GG_{\om_{a_1}}\ast\cddd\ast \GG_{\om_{a_n}}$
(here $\om_{k}$ is the $k$-th fundamental coweight of $GL(m)$)
whose $\C$-points are
$n$-step flags of lattices in $V(O)$:
$$
\TGG^a_{\mu} =\{L_0= \ L_0\sub\ L_1\sub\cddd\sub L_n=\ L ~
| ~  L\in \GG_{\mu},
\dim L_i/L_{i-1}=a_i, z(L_i)\sub L_{i-1} \},
$$
where $L_0=V(O)$.
It is known that $\pi^a_{\mu}=\pi:\TGG^a_{\mu}\ra\barr\GG_{\mu}$ is
a resolution of singularities \cite{MV}.

Consider $L\in\GG_{\la}\sub \barr\GG_{\mu}$.
Observe that
$(\pi^a_{\mu})^{-1}(L)=\FF^{n,a}_x$, where $\FF^{n,a}_x$
is the Springer-Ginzburg fiber defined in
\ref{ginzburgfiber}.

\sus{Perverse sheaves on affine Grassmannians}\label{satake}
In this subsection $G$ denotes $GL(m)$ or $PGL(m)$.

Let $\Perv_{G(O)}(\GG_{G})$
be the category of
$G(O)$-equivariant
perverse sheaves on $\GG_{G}$.
We will denote by $\IC_{\mu}=\IC(\barr\GG_{\mu})$
the intersection cohomology complex on the closure of the
orbit $\GG_{\mu}$.
There is a tensor product (convolution) construction \cite{MV, MV2}
which makes
the category $\Perv_{G(O)}(\GG_{G})$
into a tensor category.

\sss*{\bf{Theorem: geometric Satake correspondence \cite[7.1]{MV}}}
The semisimple tensor category $\Perv_{G}(\GG)$
is equivalent to the category
$\Rep G^L$
of rational representations of the
Langlands dual group $G^L$.
Under this equivalence the sheaf $\IC_{\mu}$ corresponds to the
highest weight representation $V_\mu$ of $G^L$.

Under the equivalence above the convolution
$\IC_{a_1}\ast\dots\ast\IC_{a_n}$ corresponds to the tensor product
$V_{a_1}\otimes\dots\otimes V_{a_n}$. By Gabber's decomposition theorem,
\begin{equation}\label{icdecomposition}
\IC_{a_1}\ast\dots\ast\IC_{a_n}=
\bigoplus_{\la}L_{\la}\otimes \IC_{\la}.
\end{equation}

On the level of representation theory, we have a decomposition
\begin{equation}\label{repdecomposition}
V_{a_1}\otimes\dots\otimes V_{a_n}=\bigoplus_{\la}\Mlt_{\la}\otimes V_{\la},
\end{equation}
where the sum is over all partitions $\la\leq\mu$, and $\Mlt_{\la}$
are the multiplicity vector spaces.

Taking the hypercohomology in the left and right hand side
of the equation (\ref{icdecomposition}) and comparing it to the
equation (\ref{repdecomposition}) we see that
\begin{equation}\label{multginzburg}
\Hom_{GL_{m}}(V_{a_1}\otimes\dots\otimes V_{a_n}, V_{\la})=\Mlt_{\la}
=H(\pi^{-1}(L_{\la})).
\end{equation}

\sus{Transverse slices arising from affine Grassmannians}

\sss{}\label{localtransetup}
Let us recall the setup of \ref{nilpnormalslice}: $x$
is a nilpotent operator of type $\la$ in $\End(D)$, $\dim D=N$.
Let $b=(b_1,\dots,b_m)$ be a
permutation of $\la$ (notice that $b_i\geq 1$).
Consider $b$ as a coweight of $GL(m)$
and consider the lattice
$L_b$ generated by the elements $z^{-b_i}\bfe_i$, $1\leq i\leq m$.
Clearly, $L_b\in\GG_{\la}$.

Let $D=L_b/L_0$. Then $\dim D=N$.
Define
$$
D_j=\text{span}\{e_i ~|~ b_i=j\},\qquad \text{ and }\qquad
d_j=\dim D_j.
$$
We have a decomposition of $D$ as follows:
\begin{equation}\label{decompz}
D=\bigoplus_{1\leq k\leq j\leq n-1} z^{-k}D_j.
\end{equation}

\sss{} Let us consider the group ind-scheme
$G(\C[z^{-1}])$, and let
$L^{<0}G(K)$ be subgroup of $G(\C[z^{-1}])$ which is the kernel
of the map $G(\C[z^{-1}])\to G$ defined by $z^{-1}\mapsto 0$.
Denote the $L^{<0}G(K)$-orbit of the lattice
$L_b$ in $\GG_{G}$ by $T_b$.

\sss{} We can choose a complement $L^-_{b}$ to $L_b$ such that
$V(K)=L_{b}\oplus L^-_{b}$. We define $L^-_{b}$ as the
subspace of $V(K)$ spanned by
$z^{-j}e_i$, $j> b_i$.
Denote the projection of $V(K)$
to $L_{b}$ along $L^-_{b}$ by $\pi_{b}$.

We can describe an open neighborhood $\UU^N_b$
of $L_b$ in $\GG_N$ as follows:
$$
\UU^N_b=\{ L\in\GG_N\ |\ \text{the projection } \pi_{b}:L\to L_{b}
\text{ is an isomorphism } \}.
$$

\sss{}\label{FirstElement}
We can describe the set $\UU^N_b$ in terms of certain maps,
generalizing a construction of \cite{L81}.

Any lattice $L\in \UU^N_b$ is of the form $(1+f)L_b$
where $f: L_b\to L^-_b$ is a linear map such that
$L_0\sub \ker f$. We can decompose $f$ as follows.
Let us consider the $m$-dimensional
vector space $V_b=\{z^{-b_i}\bfe_i~|~ 1\leq i\leq m\}$.
Then
$$
f=\sum_{k=1}^{\infty}z^{-k}f_k ,
$$
where $f_k:L_b/L_0\to V_b$ are linear maps.
It is easy to see that since $(1+f)L_b$ is a lattice,
we have $f_k=f_1(z+f_1)^{k-1}$ and the operator
$z+f_1:L_b/L_0\to L_b/L_0$ is nilpotent. Altogether:
\begin{equation}\label{fandf1}
f=\sum_{k=1}^{\infty}z^{-k}f_1(z+f_1)^{k-1}.
\end{equation}

Observe that if $L=(1+f)L_b$ then the isomorphism
$\pi_b$ intertwines the action of $z$ on $L/L_0$ with
the action of $z+f_1$ on $L_b/L_0$.

\sss{} Now we consider the action on $\GG_G$ of the
group of ``loop rotations'' isomorphic to the multiplicative group $\C^*$:
$z\mapsto sz,\ \  s\in \C^*$, which acts on
$V((z)) = V(K)$ by sending $z^k\bfe_i$ to $(sz)^k\bfe_i$.
Denote by $s\circ L$ the result of this action of $s\in \C^*$
on a lattice $L\in \GG_G$.

Consider this action on the lattices in $\UU^N_b$ i.e.,
lattices of the form $L=(1+f)L_b$. Our $\C^*$-action on $V(K)$ restricts
to the action on $L_b$, $L^-_b$, $L_b/L_0$ and $V_b$ and we denote $s\circ f = s\cdot f\cdot s^{-1}$
and $s\circ f_1 = s\cdot f\cdot s^{-1}$.

We have:
$$
s\circ L = s\circ (1+f)L_b=(1+s\circ f)s\circ L_b=(1+s\circ f) L_b
$$
since $L_b$ is a $T$-invariant point in $\GG_G$. Now,
$$
s\circ f=\sum_{k=1}^{\infty}(sz)^{-k} (s\circ f_1)(sz + (s\circ f_1))^{k-1}
=\sum_{k=1}^{\infty} z^{-k} s^{-1}(s\circ f_1)(z+s^{-1}(s\circ f_1))^{k-1},
$$
where $s^{-1}(s\circ f_1)$ is the composition of $(s\circ f_1)$ and the operator $s^{-1}\Id_{V_b}$ on $V_b$.

\sss{} Now the following lemma is clear:

\slemm{} $\lim_{s\to\infty}s\circ f=0$ if and only if
$\lim_{s\to\infty}s^{-1}(s\circ f_1)=0$.

\sss{} Let us now study $s\circ f_1$.
We will consider
here $f_1$ as a map from $L_b/L_0$ to itself equipped with
the basis $\{z^{-k_i}\bfe_i\ | 1\leq i\leq m, \ 1\leq k_i\leq b_i \}$.
If $u\in L_b/L_0$ is a vector then
$$
u=\sum_{k,i}u_{k,i}z^{-k}\bfe_i.
$$
Denote the matrix elements of $f_1$ is this basis by $f^{l,j}_{k,i}$
where $f^{l,j}_{k,i}:\C z^{-l}\bfe_j\to\C z^{-k}\bfe_i$.
Now recall that by construction we have
$$
f^{l,j}_{k,i}=0, \ \ \text{ if } k\neq \la_i.
$$
Then
\begin{equation}
\begin{split}
f(u)_{k,i}&= 0, \text{ if } k\neq \la_i, \\
f(u)_{\la_i,i}&= \sum_{l,j}f^{l,j}_{k,i}u_{l,j}.
\end{split}
\end{equation}
Now for $s \circ f_1=s\cdot f_1\cdot s^{-1}$ we have:
\begin{equation}\label{scircf1}
\begin{split}
(s \circ f_1)(u)_{k,i}&= 0, \ \ \text{ if } k\neq c_i, \\
(s \circ f_1)(u)_{\la_i,i}&=\sum_{l,j} s^{-\la_i}f^{l,j}_
{\la_i,i}u_{l,j} s^{l}
=\sum_{l,j}s^{l-\la_i} f^{l,j}_{\la_i,i}u_{l,j}.
\end{split}
\end{equation}

\lemm{}\label{limtslice} The following are equivalent:
\begin{enumerate}
\item $\lim_{s\to\infty}s^{-1} (s \circ f_1)=0$.
\item $f^{l,j}_{\la_i,i}=0$ if $l>\la_i$.
\end{enumerate}

\begin{proof} Follows immediately from (\ref{scircf1}).
\end{proof}

\lemm{}\label{tcstar}
A lattice $L\in \GG_G$ is in the $L^{<0}G(K)$-orbit of
$L_b$ if and only if \linebreak
$\lim_{s\to\infty} s\circ L=L_b$.

\begin{proof} We have the following decomposition \cite[Corollary 2.2]{F}:
$$
G(K)=G([z^{-1}])X_{*}(T)G(O).
$$
Then
$$
\GG=G(K)/G(O)=\bigcup_{\la\in X_{*}(T)}G([z^{-1}]) (\la\cdot G(O)).
$$
Now $G([z^{-1}])=L^{<0}G(K)G$ is a semidirect product. So,
$$
\GG=\bigcup_{\la\in X_{*}(T)}L^{<0}G(K) G (\la\cdot G(O)).
$$
The orbits of $L^{<0}G(K)$ intersect the orbits of
$G(O)$ transversally, \cite[Section 2. Remark]{F}.
This means in particular that if $p\in G\cdot\la$, then
$(L^{<0}G(K)\cdot p)\cap G\cdot\la=p$. Then we have
\begin{equation}\label{disjoint}
\GG=\bigsqcup_{\substack{\la\in X^{+}_{*}(T)\\
p\in G\cdot\la}}
L^{<0}G(K)\cdot p,
\end{equation}
where $X^{+}_{*}(T)$ is the set of dominant coweights of $G$.

Since for $g\in L^{<0}G(K)$ we have $\lim_{s\to\infty} (s\circ g)=1$,
it is clear that
$$
L^{<0}G(K)\cdot p\sub \{L\in \GG\ |\ \lim_{s\to\infty} s\circ L=p \}.
$$
Since we have the disjoint decomposition (\ref{disjoint}),
we actually have
$$
L^{<0}G(K)\cdot p= \{L\in \GG\ |\ \lim_{s\to\infty} s\circ L=p \}.
$$
\end{proof}

\lemm{}\label{sliceinopen} If $T_b:=L^{<0}G(K)\cdot L_b$ then
$T_b\cap \GG_N\sub \UU_b^N$.

\begin{proof} Clear.
\end{proof}

\sss{} Again, recall the setup of \ref{nilpnormalslice}
and the definition of the variety $T_x$.
Let $x+f_1\in T_x$. Construct a map
\begin{equation}
\begin{aligned}
\psi: T_x\cap\NN & \to \UU^N_b, \\
\psi:x+f_1 & \mapsto (1+\sum_{k=1}^{\infty}z^{-k}f_1(z+f_1)^{k-1})L_b.
\end{aligned}
\end{equation}

\lemm{} The image of $\psi$ defined above
is contained in $T_b\cap \GG_N$. Moreover,
the map $\psi: T_x\cap\NN \overset{\simeq}\ra T_b\cap \GG_N$
is an isomorphism of algebraic varieties.

\begin{proof} By definition of $T_x$,
Lemma \ref{limtslice}, and Lemma \ref{tcstar}
$$
\psi(T_x\cap\NN)=\{L\in \UU_b^N \ |\  \lim_{s\to\infty} s\circ L=L_b\}
=T_b\cap \UU_b^N.
$$
Since by Lemma \ref{sliceinopen} $T_b\cap \GG_N\sub \UU_b^N$, and
$\UU_b^N\sub \GG_N$
we have $T_b\cap \UU_b^N=T_b\cap \GG_N$.
\end{proof}

\sss{} Recall the setup of \ref{tildeta}.
Also, let $\pi=\pi^a_{\mu}:  \TGG^a_{\mu}\ra\barr\GG_{\mu}$ is a resolution
of singularities, cf. \ref{convolutionlocal}.

We can lift the map $\psi$
to the map
$$
\tpsi: \TT_{x}^{a}\ra  \pi^{-1}(T_b\cap\barr\GG_{\mu})\sub \TGG^a_{\mu}
$$
since a $(x+f_1)$-invariant $n$-step flag in $D$ will give rise
to a $n$-step flag of lattices in $L=\psi(x+f_1)$.

\slemm The map $\tpsi$ is an isomorphism of algebraic varieties.
Moreover, the following diagram of morphisms
\begin{equation}
\begin{CD}
\TT_{x}^{a}@>{\tpsi}>> \pi^{-1}(T_b\cap\barr\GG_{\mu})  \\
  @V{\mmm_a}VV  @V{\pi}VV  \\
T_{x,a} @>{\psi}>> T_b\cap\barr\GG_{\mu}
\end{CD}
\end{equation}
commutes.

\sus{Global version of the map $\psi$}\label{globalpsi}

\sss{}
Recall the setup of \ref{localtransetup}.
Let us consider the scheme $\BDG_N$
and let $\LL_b$ be the lattice in $V(z)$
generated by the elements $z^{-b_i}\bfe_i$, $1\leq i\leq m$.

Just as in the local case consider $m$-dimensional
vector subspace $V_b=\{z^{-b_i}\bfe_i~|~ 1\leq i\leq m\}$.
of $V(z)$, and consider a linear map $f_1: \LL_b/\LL_0\to V_b$.
(Notice that $D = \LL_b/\LL_0 \simeq L_b/L_0$ where $L_b$ and $L_0$
are the analogous local lattices.)

\lemm For any $u\in \LL_b/\LL_0$ and any $e\in \A^1$, we have
$$
(1+\sum_{k=1}^{\infty}z^{-k}f_1(z+f_1)^{k-1})(u)
=(1+\sum_{k=1}^{\infty}(z-e)^{-k}f_1(z-e+f_1)^{k-1})(u).
$$

\begin{proof} Binomial formula.
\end{proof}

\sss{} Now consider $z+f_1$ as an operator on $D=\LL_b/\LL_0$,
let $E$ be its spectrum and let $\pr_e: D\to D_e$, for $e\in E$,
be the projection to the generalized $e$-eigenspace.

Once again, recall the setup of \ref{nilpnormalslice}
and the definition of the variety $T_x$.
For $x+f_1\in T_x\sub \End(D)$ define the subspace $\psi(x+f_1)$ in $V(z)$
as follows
\begin{equation}\label{globallattice}
\psi(x+f_1)=(\sum_{e\in E}
(1+\sum_{k=1}^{\infty}(z-e)^{-k}f_1(z-e+f_1)^{k-1})
\pr_e) L_b.
\end{equation}

\slemm{} The subspace $\psi(x+f_1)$ is a lattice in $V(z)$ and therefore
an element in $\BDG_N$.

\begin{proof} The same as in the local case.
\end{proof}

Summarizing, we have constructed an embedding
$$
\psi: T_x\hookrightarrow \BDG_N
$$

As in the local case, this embedding lifts to an embedding
$\tpsi: \tmm^{-1}(T_x)\cap \tfg^{n,a}\hookrightarrow \TBDG^{n,a}_N$
in such a way that the diagram

\begin{equation}
\begin{CD}
\tmm^{-1}(T_x)\cap \tfg^{n,a}  @>{\tpsi}>{\subset}> \TBDG^{n,a}_N  \\
  @V{\mmm_a}VV  @V{\pi}VV  \\
T_{x} @>{\psi}>{\subset}>  \BDG_N
\end{CD}
\end{equation}
commutes.

\se{Main Results}\label{mainresults}

\sus{Combinatorial data}\label{data}

\sss{From quiver data to $GL(n)$-data}\label{qdatatogl}
Let $d=(d_1,\dots, d_{n-1})$ and $v=(v_1,\dots, v_{n-1})$
be two $(n-1)$-tuples of non-negative integers. We will
transform this "quiver data" into some $GL(n)$ weights.
\begin{enumerate}
\item Let $C$ be the Cartan matrix of type $A_{n-1}$. By $(d - Cv)_j$
we will denote the $j$-th component of the $(n-1)$-tuple $d-Cv$.
\item Let $N=\sum_{j=1}^{n-1}jd_j$ and let $m=\sum_{j=1}^{n-1}d_j$.
\item Let $\cla=(\cla_1, \cla_2,\dots,\cla_n)$ be a partition
of $N$ defined as follows (here $d_n=0$):
$$
\cla_i=\sum_{j=i}^{n}d_j.
$$
\item Let $\la$ be the dual partition.
\item Let $a=(a_1,\dots, a_n)$ be defined as follows,
cf. \cite[8.3]{N94}, (here $(d-Cv)_n=0$):
\begin{equation}\label{a}
a_i=v_{n-1}+\sum_{j=i}^{n}(d-Cv)_j.
\end{equation}
\item Let $\cmu$ be the partition obtained from $a$ by
permutation and let $\mu$ be the dual partition.
\end{enumerate}
We can view $\cla$ as a highest weight of $GL(n)$ and $a$
as a weight in the highest weight $GL(n)$-module $W_{\cla}$,
cf. \ref{fromSLntoGLn}.

\sss{From $GL(n)$-data to a conjugacy class}\label{gldatatoclass}

Let $c=(c_1,\dots,c_{n-1})$ be in the center of $\fg(V)$, where
$\fg(V)=\prod_{i=1}^{n-1}\gl(V_i)$ and $\dim V_i=v_i$.

First of all, denote
\begin{equation}
b_1=0\ \text{ and }\ b_i=c_1+\dots +c_{i-1}\ \text{ for }\ 2\leq i\leq n.
\end{equation}
Let $b=(b_1,\dots, b_{n})\subset \A^n$. Let
$P$ be the polynomial $P(t)=\prod_{i=1}^n(t-b_i)$.

Consider $E=E(c)=\{b_1,\dots, b_{n}\}$ as a subset of $\A^1$
and consider $b$ as a map $[1,n]\to E$ defined by $b(i)=b_i$.

For every $e\in E$ denote
$I(e):= b^{-1}(e)=\{i\in[1,n]~|~ b_i=e\}$.
Now take $a$ as in (\ref{a}) and let $a(e)=(a_i)_{i\in I(e)}$
and let $\cmu(e)$ be the partition obtained from $a(e)$
by permutation. Let $\mu(e)$ be the dual partition,
and let $\tmu=\{\mu(e)\}_{e\in E}$ be the collection
of all partitions attached to eigenvalues.
Let $\OO_{E,\tmu}$ be the conjugacy class in $\End(D)$, $\dim D=N$
attached to the data $E,\tmu$ as in \ref{conjclass}.

\sus{} Now we can formulate our main theorem. For notation
on quiver varieties see \ref{quiverdefs},
on Springer-Ginzburg resolutions
see \ref{definitionBases},
on transverse slices see
\ref{sliceJordanBasis}, and finally on
Beilinson-Drinfeld Grassmannians see \ref{globalpicture}.

\Theo{}\label{main} Let $N,m,v,d,a,c,b,E,\la,\tmu$ be as above.
There exist algebraic isomorphisms
$\phi, \tphi$ and algebraic immersions  $\psi, \tpsi$
such that the following diagram commutes:
\begin{equation}
\begin{CD}
\fM(v,d)
@>{\tphi}>{\simeq}>
\tmm^{-1}(T_{\la})\cap \tfg^{n,a,E,\tmu}
@>{\tpsi}>{\subset}> \TBDG^{n,a}_b(P)
\\
@V{p}VV
@V{\tmm}VV
@V{\pi}VV
\\
\fM_1(v,d)
@>{\phi}>{\simeq}>
T_{\la}\cap \barr\OO_{E,\tmu}
@>{\psi}>{\subset}>
\BDG_{N,b}(P).
\end{CD}
\end{equation}

\sus{Remarks and Corollaries}

\sss{Remark} When $c=0$ we can describe the images of the maps
$\psi$ and $\tpsi$ and obtain a more precise result
stated in the introduction and \cite{MVyb}. In particular,
$(\psi\circ\phi)(0)=L_\la\in \GG_0$,
and
$\tpsi\circ\tphi$
restricts to an isomorphism
\begin{equation}\label{fiberiso}
\tpsi\circ\tphi: \fL(v,d)\simeq\pi^{-1}(L_\la).
\end{equation}
We believe that one should be able to generalize these statements for arbitrary $c$.

\sss{Dimensions} Let $c=0$.
First of all we'll check that the varieties
$\fM(v,d)$ and $\TT_{x}^{a}$ have the same dimension.
According to Nakajima \cite[Corollary 3.12]{N98}
$\fM(v,d)$, if nonempty, is a smooth variety of dimension
${}^{t}v(2d-Cv)$
where $C$ is the Cartan matrix of type $A_{n-1}$.
If $\cla$ and $\cmu$ are defined by $v,d$ as in \ref{qdatatogl},
then we have
\begin{equation}
\begin{split}
\dim\fM(v,d) & ={}^{t}v(2d-Cv)=
2\sum_{i=1}^{n-1}v_id_i-2\sum_{i=1}^{n-1}v_i^2
+2\sum_{i=1}^{n-2}v_iv_{i+1}  \\
& =\sum_{i=1}^{n-1}[(\cla_i)^2-(\cmu_i)^2]=
\dim\TT_{x,\mu}.
\end{split}
\end{equation}

We will list here two applications of our Main Theorem.

\sss{\bf A compactification of quiver varieties}
The closure in $\BDG_{N,b}(P)$
of the image of $\fM_1(v,d)$ under the map $\psi\circ\phi$
gives us a compactification of $\fM_1(v,d)$.
Analogously, the closure in $\TBDG^{n,a}_b(P)$
of the image of $\fM(v,d)$ under the map $\tpsi\circ\tphi$
gives us a compactification of the quiver variety
$\fM(v,d)$.

\sss{\bf A decomposition of the affine Grassmannian}\label{decompositionAffGrass}

The following is a corollary of the main theorem. Here $c=0$.

\scorr{} We can decompose $\barr\GG_{\mu}$ into the following disjoint union:
\begin{equation}
\barr\GG_{\mu}=\bigsqcup_
{\substack{y\in G\cdot\la \\ \la\leq\mu  }}
\fM_0(v,d)_y ,
\end{equation}
where $\la$ varies over the set of dominant coweights of $G$,
$G\cdot\la$ is the $G$-orbit of $\la$ in $\GG_G$, and
$\fM_0(v,d)_y$ is a copy of quiver variety $\fM_0(v,d)$
for every point $y\in G\cdot\la$,
with $v,d$ obtained from $\la,\mu$
by reversing the procedures of \ref{qdatatogl}.

\begin{proof} As in the proof of \ref{tcstar}, we have:
\begin{equation}
\GG_G=\bigsqcup_{\substack{\la\in X^{+}_{*}(T)\\
y\in G\cdot\la}}
L^{<0}G(K)\cdot y .
\end{equation}
Then:
\begin{equation}
\barr\GG_{\mu}=\bigsqcup_{\substack{\la\in X^{+}_{*}(T)\\
y\in G\cdot\la}}
(L^{<0}G(K)\cdot y)\cap\barr\GG_{\mu}=
\bigsqcup_{\substack{\la\in X^{+}_{*}(T)\\
y\in G\cdot\la}} \fM_0(v,d)_y
\end{equation}
since every $(L^{<0}G(K)\cdot y)\cap\barr\GG_{\mu}$, for $y\in G\cdot\la$
is isomorphic to a copy of $\fM_0(v,d)$.
\end{proof}

\sss{Remarks}
\begin{enumerate}
\item An ``affine analogue" of our construction has recently appeared in the paper \cite{BF}.
\item We would also like to mention another example of a decomposition
of an infinite Grassmannian into a disjoint union of quiver varieties.
Generalizing a result of G. ~Wilson \cite{W},
V. ~Baranovsky, V.~ Ginzburg,  and A.~ Kuznetsov \cite{BGK}
constructed a decomposition of (a part of) adelic Grassmannian
into a disjoint union of \emph{deformed} versions of quiver varieties
$\fM(v,d)$ associated to affine quivers of type A.
\end{enumerate}

\se{On quiver varieties and conjugacy classes of matrices}\label{quiverConjClassSection}

\sus{Definitions} Let us consider a particular case
of the Main Theorem. Let $d=(N, 0,\dots, 0)$
and $v=(v_1,\dots, v_{n-1})$
be the $(n-1)$-tuple of non-negative integers
such that $N\geq v_1\geq v_2\geq\dots\geq v_{n-1}$.

\sss{}
Define the algebraic morphisms $\tphi: \fM(v,d)\ra\tfg^{n,a,E,\tmu}$ and
$\phi: \fM_1(v,d)\ra\barr\OO_{E,\tmu}$
as follows:
\begin{equation}
\begin{aligned}
\tphi: (x,\bx,p,q)& \mm (q_1p_1, \{0\}\sub\ker p_1\sub\ker x_1p_1
\sub\ker x_{n-1}\dots x_1p_1) , \\
\phi:(x,\bx,p,q)& \mm q_1p_1 .
\end{aligned}
\end{equation}

The following theorem is a common generalization of
(some of) the results of \cite{KP} and \cite{N94},
cf. \cite{CB}.

\Theo{}\label{gkp} The maps $\phi,\tphi$ defined above are
isomorphisms of algebraic varieties and the following diagram commutes

\begin{equation}
\begin{CD}
\fM(v,d) @>{\tphi}>> \tfg^{n,a,E,\tmu}\\
@V{p}VV @V{\tmm}VV \\
\fM_1(v,d) @>{\phi}>> \barr\OO_{E,\tmu}
\end{CD}
\end{equation}

\begin{proof}  Following the logic of \cite{N98, M},
it is not hard to check that $\tphi$ is a bijective
morphism between two smooth varieties of the same dimension
and thus an isomorphism. The map $\phi$ is a closed immersion
and it is surjective since both $p$ and $\tmm$ are surjective.
\end{proof}

\sss{}
In particular, if all the numbers
$0, c_1, c_1+c_2, \dots, c_1+c_2+\dots+ c_{n-1}$
are pairwise distinct, then the quiver variety $\fM(v,d)$
is isomorphic to the conjugacy class of a semisimple element
(diagonal matrix)
$$
\operatorname{diag}
(b_1,\dots, b_1,b_2,\dots, b_2,  \dots,
b_{n},\dots, b_{n}),
$$
where $b_1=0$ appears with multiplicity $a_1$,
$b_2=c_1$ appears with multiplicity $a_2$, and so on,
and $b_{n}=c_1+c_2+\dots, c_{n-1}$ appears with multiplicity $a_n$.

\sss{Remark} In fact one can also prove that the quiver variety
$\fM_0(v,d)$ is isomorphic to a conjugacy class which
is generally different from the conjugacy class considered above.
The two classes coincide when the $SL(n)$ weight $d-Cv$ is dominant,
i.e. when $a_1\geq a_2\geq \dots \geq a_n$.

\se{Proof of the Main Lemma}\label{proofMainLemma}

\sus{D'apr\` es Maffei}\label{dapresMaffei}

\sss{}\label{mafNotation} We borrow Maffei's \cite{M} notations and conventions.
Let $v=(v_1,\dots,v_{n-1})$ and
$d=(d_1,\dots,d_{n-1})$ be two $(n-1)$-tuples of integers
and let us define $(n-1)$-tuples $\tv$ and $\td$ as follows:
\begin{equation}
\begin{split}
\td_1 &:=\sum_{j=1}^{n-1}jd_j, \\
\td_i &:=0, \text{ for } i>1, \\
\tv_i &:=v_i+\sum_{j=i+1}^{n-1}(j-i)d_j. \\
\end{split}
\end{equation}

Our goal is to construct a map from $\La^c(v,d)$ to $\La^c(\tv,\td)$,
that is we have to send a quadruple $(x,\bx,p,q)\in \La^c(v,d)$
to a quadruple $(\TA,\TB,\tga,\tde)\in \La^c(\tv,\td)$.
First of all,
the $I$-graded vector spaces $\TV_i$ and $\TD_i$ such that
$\dim\TV_i=\tv_i$ and $\TD_i=\td_i$
are constructed as follows. Let $D_j^{(k)}$ be a copy of $D_j$.
\begin{equation}
\begin{split}
\TD_1 &=\bigoplus_{1\leq k\leq j\leq n-1} D_j^{(k)}, \\
\TD_i &=0, \text{ for } i>1, \\
\TV_i &=V_i\oplus\bigoplus_{1\leq k\leq j-i\leq n-i-1} D_j^{(k)}.\\
\end{split}
\end{equation}

We need the following subspaces of $\TV_i$.

\begin{equation}
D _i ^{\prime} =
\bigoplus_{\substack{  i+1 \leq j \leq n-1 \\
1 \leq k \leq   j-i}}  D_j^{(k)},  \qquad
D _i ^{+} = \bigoplus_{\substack{i+2 \leq j \leq n-1 \\
2 \leq k \leq   j-i}} D_j^{(k)},   \qquad
D _i ^{-} = \bigoplus_{\substack{ i+2 \leq j \leq n-1 \\
1 \leq k \leq   j-i-1}}  D_j^{(k)}.
\end{equation}

In order to make the notation more homogeneous we set
$\TV_0:=\TD_1$, $\TA_0=\tga_1$, $\TB_0=\tde_1$.

We will name the blocks of the maps $\TA_i$ and $\TB_i$ as follows
\begin{equation}\label{blocks}
\begin{aligned}
\pi_{D^{(h)}_j} \TA_i|_{D^{(h')}_{j'}}
&= {}^{i}\!t^{j',h'}_{j,h} &
\qquad   \pi_{D^{(h)}_j} \TB_i|_{D^{(h')}_{j'}}
&= {}^{i}\! s^{j',h'}_{j,h}  \\
\pi_{D^{(h)}_j} \TA_i|_{V_i } &= {}^{i}\! t^{V}_{j,h} &
\qquad \pi_{D^{(h)}_j} \TB_i|_{V_{i+1} } &= {}^{i}\!  s^{V}_{j,h} \\
\pi_{V_{i+1}}\TA_i|_{D^{(h')}_{j'}} &= {}^{i}\!  t^{j',h' }_{V} &
\qquad
\pi_{V_{i}}\TB_i|_{D^{(h')}_{j'}} &= {}^{i}\!s^{j',h'}_{V}
\end{aligned}
\end{equation}

We define also the following operator $z_i$ on $D_i^{\prime}$
\begin{equation}
\begin{aligned}
z_i|_{D_j^{(1)}} &= 0   ,                         \\
z_i|_{D_j^{(h)}} &= Id_{D_j}: D_j^{(h)} \ra D_j^{(h-1)}
\end{aligned}
\end{equation}

\sss{} Following Maffei let us introduce the following degrees:
\begin{equation}
\begin{aligned}
 \deg({}^{i}\!t^{j',h'}_{j,h} )  &=  \min(h-h'+1, h-h'+1+j'-j), \\
 \deg({}^{i}\!s^{j',h'}_{j,h} )  &=  \min(h-h', h-h'+ j'-j).
\end{aligned}
\end{equation}

\sss{} A quadruple $(\TA,\TB,\tga,\tde)\in \La^c(\tv,\td)$
is called \emph{transversal} if it satisfies the following two groups
of relations for $0\leq i\leq n-2$
\begin{enumerate}
\item first group (Maffei)
\begin{equation} \label{transone}
\begin{aligned}
   {}^{i}\!t^{j',h'}_{j,h}  &= 0
& &\mbox{ if } \deg(t^{j',h'}_{j,h}) < 0  \\
   {}^{i}\!t^{j',h'}_{j,h}  &= 0
& &\mbox{ if } \deg(t^{j',h'}_{j,h}) = 0
\mbox{ and } (j',h') \neq (j, h+1) \\
   {}^{i}\!t^{j',h'}_{j,h}  &= Id_{D_j}
& &\mbox{ if } \deg(t^{j',h'}_{j,h}) = 0
\mbox{ and } (j',h') = (j, h+1) \\
   {}^{i}\!t^{V}_{i,j,h}        &= 0         & &   \\
   {}^{i}\!t^{j',h'}_{V}    &= 0
& &\mbox{ if } h' \neq 1 \\
   {}^{i}\!s^{j',h'}_{j,h}  &= 0
& &\mbox{ if } \deg(s^{j',h'}_{j,h}) < 0  \\
   {}^{i}\!s^{j',h'}_{j,h}  &= 0
& &\mbox{ if } \deg(s^{j',h'}_{j,h}) = 0
\mbox{ and } (j',h') \neq (j, h) \\
   {}^{i}\!s^{j',h'}_{j,h}  &= Id_{D_j}
& &\mbox{ if } \deg(s^{j',h'}_{j,h}) = 0
\mbox{ and } (j',h') = (j, h) \\
   {}^{i}\!s^{V}_{j,h}        &= 0
& &\mbox{ if } h \neq j-i  \\
   {}^{i}\!s^{j',h'}_{V}    &= 0         & &
\end{aligned}
\end{equation}
\item second group
\begin{equation}\label{transtwo}
\pi_{D^{(h)}_j} \TB_i \TA_i|_{D^{(h')}_{j'}} - x_i=0
\qquad \text{ unless } h=j-i
\end{equation}
\end{enumerate}

Let us denote the set of all transversal elements in
$\La^c(\tv,\td)$ by $S$.
The set of all stable transversal elements is denoted by
$S^s=S\cap\La^{c,s}(\tv,\td)$.

\sss{} We will need more notation. First of all denote
\begin{equation}
b^{i}_{j}=c_{i+2}+\dots + c_{j}\qquad \text{ for } -1\leq i\leq n-3,
\text{ and } i+2\leq j\leq n-1.
\end{equation}

Now we introduce some invariant polynomials of
$q_{i\ra j}p_{j\ra i}$ as follows.
First,
\begin{equation}
P(i,1,j)=q_{i+2\ra j}p_{j\ra i+2}
\end{equation}
and for $2\leq h'\leq j-i-1$
\begin{equation}
\begin{split}
P(i,h',j) & =q_{i+h'+1\ra j}p_{j\ra i+h'+1} \\
 & +\sum_{k=1}^{j-i-h'-1} (-1)^k
\si_k(b^{i}_{i+2},\dots,b^{i}_{i+h'-1+k})
q_{i+h'+1+k\ra j}p_{j\ra i+h'+1+k} \\
 & +(-1)^{j-i-h'-1}
\si_{j-i-h'}(b^{i}_{i+2},\dots,b^{i}_{j-1}).
\end{split}
\end{equation}
where $\si_k$ is the $k$-th elementary symmetric function.

We also fix the notation for binomial coefficients
$$
\left(\begin{matrix}
n \\ k
\end{matrix}\right)=\frac{n!}{k!(n-k)!}.
$$

\sus{Main Lemma} We can now formulate our main lemma

\slemm\begin{enumerate}
\item[(i)] There exists a unique $G(V)$-equivariant map
$\Phi: \La^c(v,d)\ra S$
such that
\begin{equation}\label{basic}
\begin{aligned}
\pi_{V_{i+1}} \TA_i|_{V_i} &= x_i         &
\pi_{V_{i}} \TB_i|_{V_{i+1}}  &= \bx_i   \\
{}^{i}\!t^{i+1,1}_{V}      &=  p_{i+1} & \ \ \ \ \
{}^{i}\!s^{V}_{i+1,1}      &=  q_{i+1}
\end{aligned}
\end{equation}
\item[(ii)] The blocks of $\TA_i, \TB_i$ not defined in the equations
(\ref{transone}) and (\ref{basic})
are described as follows:
\begin{equation}
\begin{aligned}
{}^{i}\!t^{j',1}_{V}      &=  p_{j'\ra i+1} &\ \ \ \ \
{}^{i}\!s^{V}_{j,j-i}      &=  q_{i+1\ra j} \label{V} \\
\end{aligned}
\end{equation}
When $j'\neq j$ we have
\begin{equation}\label{tj}
\begin{aligned}
{}^{i}\!t^{j',h'}_{j,h} & =0  &
\text{ if } & (j',h') \neq (j, h+1) \\
{}^{i}\!s^{j',h'}_{j,h}  &= 0  &  \text{ if } &
(j',h') \neq (j, h) \text{ and } h\neq j-i
\end{aligned}
\end{equation}
and
\begin{equation}\label{sj}
{}^{i}\!s^{j',h'}_{j,j-i}=
q_{i+h'+1\ra j}p_{j'\ra i+h'+1}
\end{equation}
When $j=j'$ we have
\begin{equation}
{}^{i}\!t^{j,h'}_{j,h}  =
\begin{cases}
0, & \text{ if } h'=1 \\
(-1)^{h-h'+1}
\left(\begin{matrix}
h-1 \\ h'-2
\end{matrix}\right) c_{i+1}^{h-h'+1},
& \text{ if } 2\leq h'\leq h+1
\end{cases}
\end{equation}
And finally,
\begin{equation}
{}^{i}\!s^{j,h'}_{j,h}=
\begin{cases}
\left(\begin{matrix}
h-1\\ h'-1
\end{matrix}\right)
c_{i+1}^{h-h'}, &
\text{ if } h\neq j-i \\
P(i,h',j)+
\left(\begin{matrix}
h-1\\ h'-1
\end{matrix}\right)
c_{i+1}^{h-h'},
& \text{ if } 1\leq h'\leq h,\text{ and }  h=j-i
\end{cases}
\end{equation}
\item[(iii)] For $x\in\La^c(v,d)$ we have $\Phi(x)\in S^s$ if and only if
$x\in\La^c_s(v,d)$. Thus the restriction of $\Phi$ to the stable points
provides the $G(V)$-equivariant map $\Phi^s: \La^c_s(v,d)\ra S^s$
\item[(iv)] The maps $\Phi$ and $\Phi^s$ are isomorphisms of algebraic varieties.
\end{enumerate}

\begin{proof}
Following Maffei, we prove the lemma by decreasing induction on $i$.
If $i=n-2$ the maps $\TA_{n-2}$ and $\TB_{n-2}$ are completely defined by
the relations (\ref{basic}) and (\ref{transone}) and it is easy to see that
$\TA_{n-2}\TB_{n-2}=c_{n-1}$.

Assume that $\TA_k, \TB_k$ are defined for $k>i$
by the formulas in the lemma.

We have the following equations for
$\TA_i$ and $\TB_i$:
\begin{equation}\label{m=l+c}
\TA_i\TB_i =\TB_{i+1}\TA_{i+1} +c_{i+1}
\end{equation}
\begin{equation}\label{n=0}
\pi_{D^{(h)}_j} \TB_i \TA_i|_{D^{(h')}_{j'}} - z_i=0\
\text{ unless } h=j-i.
\end{equation}
Observe that
\begin{equation}\nonumber
\pi_{V_{i+1}} \TA_i \TB_i |_{V_{i+1}}= A_iB_i+p_{i+1}q_{i+1}
= B_{i+1}A_{i+1} +c_{i+1}= \pi_{V_{i+1}} \TB_{i+1} \TA_{i+1} |_{V_{i+1}}
+c_{i+1}.
\end{equation}
Then, in agreement with formulas (\ref{V})
\begin{equation}\nonumber
\begin{aligned}
\pi_{V_{i+1}} \TA_i \TB_i |_{D_{j}^{(h)}} =
\pi_{V_{i+1}} \TB_{i+1} \TA_{i+1} |_{D_j^{(h)}}=
K_{h,1} B_{i+1} p_{j\ra i+2}=
K_{h,1} p_{j\ra i+1}, \notag \\
\pi_{D_{j}^{(h)}} \TA_i \TB_i |_{V_{i+1}} =
\pi_{D_{j}^{(h)}} \TB_{i+1} \TA_{i+1} |_{V_{i+1}}=
K_{h,j-i-1} q_{i+2\ra j} A_{i+1}
=K_{h,j-i-1}q_{i+1 \ra j}.\notag
\end{aligned}
\end{equation}
where
$$
K_{p,q}=\begin{cases}
1,\  p=q \\
0,\ p\neq q
\end{cases}
$$
is the Kronecker symbol.

Now, in order to simplify the notation a bit we set
$t^{j',h'}_{j,h}:=\ {}^{i}\!t^{j',h'}_{j,h}$
and $s^{j',h'}_{j,h}:=\ {}^{i}\!s^{j',h'}_{j,h}$

Case I: $j\neq j'$. In this case
the equation (\ref{m=l+c}) and translates into the following
equations for $t^{j',h'}_{j,h}$ and $s^{j',h'}_{j,h}$:
\begin{equation}\label{m=lts}
s^{j',h'}_{j,h+1}+\sum_{\substack{h'<h''<h+1 \\ h'-j'<h''-j''<h+1-j}}
t^{j'',h''}_{j,h}s^{j',h'}_{j'',h''}+t^{j',h'}_{j,h}=
\begin{cases}
0, & \mbox { if } h\neq j-i-1 \\
q_{i+h'+1\ra j}p_{j'\ra i+h'+1}, & \mbox { if } h= j-i-1. \\
\end{cases}
\end{equation}
while the equation (\ref{n=0}) translates into the following
equations for $t^{j',h'}_{j,h}$ and $s^{j',h'}_{j,h}$, $h\neq j-i$:
\begin{equation}\label{n=0ts}
t^{j',h'}_{j,h}+\sum_{\substack{h'-1<h''<h \\
h'-1-j'<h''-j''<h-j}}s^{j'',h''}_{j,h}t^{j',h'}_{j'',h''}
+s^{j',h'-1}_{j,h}=0
\end{equation}

We claim that the system of equations (\ref{m=lts}) and (\ref{n=0ts})
has a unique solution indicated in the statement of the lemma.
We will prove this claim by
induction on $h$ and $h'$.

First of all, observe that from the
equation (\ref{n=0ts}) we have $t^{j',1}_{j,1}=0$.

We make two induction assumptions ($k\geq 1$):
\begin{enumerate}
\item  $t^{j',h'}_{j,h}=0$ for all $(h',h)$
such that $h'\leq h\leq k$
for all $j\neq j'$ at the same time.
\item  $s^{j',h'}_{j,h+1}=0$ for all $(h',h)$ such that
$h'< h\leq k+1\leq j-i$
for all $j\neq j'$ at the same time.
\end{enumerate}

Induction Step 1. Consider the equation (\ref{n=0ts}) for $h=k+1$.
By assumption (2) we have $s^{j',h'-1}_{j,k+1}=0$ and
$s^{j'',h''}_{j,k+1}=0$ for $j''\neq j$. If $j''=j$, then $j''\neq j'$
and by assumption (1) $t^{j',h'}_{j'',h''}=0$ for $h''\leq k$.
Now from equation (\ref{n=0ts}) we see that $t^{j',h'}_{j,k+1}=0$
for $h'\leq k+1$.

Induction Step 2. Consider the equation (\ref{m=lts}) for $h=k+1$.
By induction step (1) $t^{j',h'}_{j,k+1}=0$ and
$t^{j'',h''}_{j,k+1}=0$ for $j''\neq j$. If  $j''=j$, then $j''\neq j'$
and by assumption (2) $s^{j',h'}_{j'',h''}=0$.
Now from equation (\ref{m=lts}) we see that $s^{j',h'}_{j,k+2}=0$
for $h'< k+2$.

Finally, if $h+1=j-i$, then the equations (\ref{m=lts})
and the induction steps 1 and 2 yield:
\begin{equation}
s^{j',h'}_{j,j-i}=
q_{i+h'+1\ra j}p_{j'\ra i+h'+1}.
\end{equation}

Case II: $j=j'$. In this case we fix $j$ and
simplify the notation further a bit,
by setting $t^{h'}_{h}:=t^{j,h'}_{j,h}$ and
$s^{h'}_{h}:=s^{j,h'}_{j,h}$.
Now, taking into account Case I,
the equation (\ref{m=l+c}) and translates into the following
equations for $t^{h'}_{h}$ and $s^{h'}_{h}$:
\begin{equation}\label{m=ltsj}
s^{h'}_{h+1}+\sum_{h'<h''<h+1}
t^{h''}_{j,h}s^{h'}_{h''}+t^{h'}_{h}=
\begin{cases}
0, & \text { if } h\neq j-i-1 \text { and } h\neq h' \\
c_{i+1} & \text { if } h\neq j-i-1 \text { and } h= h' \\
P(i,h',j), & \mbox { if } h= j-i-1 \text { and } h\neq h'\\
P(i,h',j)+c_{i+1}, & \mbox { if } h= j-i-1 \text { and } h=h'\\
\end{cases}
\end{equation}
(In order to compute the right hand side, we need to use the
following combinatorial formula
$$
\si_a(c, c+b_1, \dots, c+b_p)
=\sum_{l=0}^a c^l
\left(\begin{matrix}
p-a+l+1 \\ l
\end{matrix}\right)
\si_{a-l}(b_1,\dots, b_p)
$$
for $a, p\in \Z$, $1\leq a\leq p$. We assume here that $\si_0(b_1,\dots, b_p) = 1$.)

The equation (\ref{n=0}) translates into the following
equations for $t^{j',h'}_{j,h}$ and $s^{j',h'}_{j,h}$, $h< j-i$:
\begin{equation}\label{n=0tsj}
t^{h'}_{h}+\sum_{h'-1<h''<h}
s^{h''}_{h}t^{h'}_{h''}
+s^{h'-1}_{h}=0
\end{equation}

Again, we claim that the system of equations (\ref{m=ltsj}) and
(\ref{n=0tsj})
has a unique solution indicated in the statement of the lemma.
Again, we will prove this claim by
induction on $h$ and $h'$.

First of all, observe that from the
equation (\ref{n=0tsj}) we have $t^{1}_{1}=0$.

We make two induction assumptions ($k\geq 1$):
\begin{enumerate}
\item $t^{h'}_{h}$ is given
by equations (\ref{tj})
 for all $(h',h)$
such that $h'\leq h\leq k$.
\item $s^{h'}_{h}$
is given
by equations (\ref{sj})
for all $(h',h)$ such that
$h'< h\leq k+1\leq j-i$.
\end{enumerate}

Proceeding by induction as in Case I
and using the formula (for $b, l\in \Z$, $0\leq b\leq l-2$)
$$
\sum_{a=b}^l (-1)^{l-a}
\left(
\begin{matrix}
l \\ a
\end{matrix}\right)
\left(\begin{matrix}
a+1 \\ b+1
\end{matrix}\right)
=0
$$
it is easy to see that all $t^{h'}_{h}$ and $s^{h'}_{h+1}$ are
given by formulas (\ref{tj}) and (\ref{sj}) respectively.

We have proved the assertions (i) and (ii)
of the lemma.
The assertion (iii) follows from the construction
and Lemma \ref{stablelemma} exactly as in \cite[Lemma 19]{M}.
The assertion (iv) follows from the construction, cf.
\cite[Lemma 19]{M}.\end{proof}

\sss{} It is important for us to record the formula
for $\TB_0\TA_0=\tde_1\tga_1$.
To simplify notation, we set
$$
b_{l}:=b^{-1}_l=c_1+\dots + c_l,
$$
and
$$
P'(h',j):=\sum_{k=1}^{j-h'-1} (-1)^k
\si_k(b_{1},\dots, b_{h'-2+k})
q_{h'+k\ra j}p_{j\ra h'+k}
+(-1)^{j-h'-1}
\si_{j-h'}(b_{1},\dots,b_{j-1}).
$$

Now we have
\begin{equation}\label{dega1c}
(\tde_1\tga_1)^{j',h'}_{j,h}=
\begin{cases}
\Id_{D_j}, &\text{ if } h'=h+1,\ j'=j, \\
q_{h'\ra j}p_{j'\ra h'}+ K_{j,j'}P'(h',j),
& \text{ if } h=j, \\
0, & \text{ otherwise }.
\end{cases}
\end{equation}
where $K_{p,q}$ is the Kronecker symbol,
$\si_k$ is the $k$-th elementary symmetric function,
and we assume that the value of $\si_k$ at the empty
collection of variables is zero.

Finally, let us record the specialization of the above
formula for the case $c=0$.
Clearly, in this case $P'(h',j)=0$ and we have
\begin{equation}\label{dega1}
(\tde_1\tga_1)^{j',h'}_{j,h}=
\begin{cases}
\Id_{D_j}, &\text{ if } h'=h+1,\ j'=j, \\
q_{h'\ra j}p_{j'\ra h'},  & \text{ if } h=j, \\
0, & \text{ otherwise. }
\end{cases}
\end{equation}

\se{Proof of the Main Theorem}\label{proofMainTheorem}

In this section we complete the proof of the Main Theorem (Theorem \ref{main}.)

\sus{The isomorphisms $\phi$ and $\tphi$}
The argument in this subsection is for the case $c=0$.
The argument for a general $c$ is completely
analogous. In the proof we mostly follow the logic of
\cite{M}.

\slemm{} Let $(\TA,\TB,\tga,\tde)\in S$ and let $\tg\in G(\TV)$
be such that $\tg(\TA,\TB,\tga,\tde)\in S$. Then $\tg_i(V_i)\sub V_i$
and if we denote $g_i=\tg_i|_{V_i}$ we have
\begin{equation}
\tg(\TA,\TB,\tga,\tde)=g(\TA,\TB,\tga,\tde) .
\end{equation}

\begin{proof} The proof is lifted verbatim from \cite[Lemma 22]{M}.
\end{proof}

\sss{} Let $D = \TD_1$ as in \ref{mafNotation}. Then $\dim D = N = \td_1 :=\sum_{j=1}^{n-1}jd_j$.
Observe that $(\tv,\td)$ as constructed in \ref{mafNotation} must satisfy the conditions of Section
\ref{quiverConjClassSection} in order for $\fM_1(\tv,\td)$ and $\fM(\tv,\td)$ to be nonempty, cf. \cite[1.4]{M}
and therefore, if nonempty,  $\fM_1(\tv,\td) \simeq \barr\OO_\mu$ and
$\fM(\tv,\td) \simeq T^*\FF^{n,a}$ (for $c = 0$), where $\mu, a$ are defined as in \ref{qdatatogl}.
(For a general $c$ the nilpotent orbit $\OO_\mu$ deforms into a general conjugacy class, cf.
\ref{conjclass}, \ref{gldatatoclass}.)
Now recall the definition (cf. \ref{nilpnormalslice}) of the transverse slice $T_x$
to the orbit $\OO_\la$ where $\la$ is obtained from $(v, d)$ as in \ref{qdatatogl}.
Let $T_{x,\mu} = T_x\cap\barr\OO_\mu$ be as in \ref{nilpnormalslice} and let $\TT_x^a$
be as in \ref{tildeta}.

\sss{}\label{constructphi} Now we will construct the maps $\phi_0$ and $\tphi$
completing the following commutative diagrams.

\begin{equation}
\begin{CD}
\La^c(v,d)@>{\Phi}>> S \\
@VVV @VVV \\
\fM_0(v,d)@>{\phi_0}>> \fM_0(\tv,\td) \\
\end{CD}
\qquad
\begin{CD}
\La^c_s(v,d)@>{\Phi^s}>> S^s \\
@VVV @VVV \\
\fM(v,d)@>{\tphi}>>  \fM(\tv,\td)\\
\end{CD}
\end{equation}

We denote $\phi:=\phi_0|_{\fM_1(v,d)}:\fM_1(v,d)\to \fM_0(\tv,\td)$.
Since $\fM_1(\tv,\td) \simeq \barr\OO_\mu$ an element of $\fM_1(v,d)$ will be sent
by $\phi$ to an operator $y + f\in \End(D)$, where y is nilpotent of type $\la$ and $f$ is given by the
\emph{explicit formulas} (\ref{dega1}) (and (\ref{dega1c}) for arbitrary $c$). A simple inspection
shows that $\Im\ \phi\sub T_{x,\mu}$, and $\Im\ \tphi\sub\TT_x^a$.

\lemm{} The map $\phi$ is a closed immersion.

\begin{proof} It is enough to prove that $\phi_0$ is closed immersion.
Recall that
\begin{equation}
\begin{split}
\fM_0(v,d) & =\La^c(v,d)//G(V)=\Spec \RR(\La^c(v,d))^{G(V)}, \\
\fM_0(\tv,\td) & =\La^c(\tv,\td)//G(\TV)=\Spec \RR(\La^c(\tv,\td))^{G(\TV)}.
\end{split}
\end{equation}
We will prove that the restriction map
$\phi^*: \RR(\La^c(\tv,\td))^{G(\TV)}\to\RR(\La^c(v,d))^{G(V)}$
is surjective.

By Theorem \ref{invtheorem} the algebra
$\RR(\La^c(\tv,\td))^{G(\TV)}$ is generated by
$\widetilde\chi(\tde_1\tga_1)$ where
$\widetilde\chi$ is a linear form
on $\Hom(\TD_1,\TD_1)$. If $\tde_1\tga_1$ is of the form
(\ref{dega1}) and
$$
\widetilde\chi=\chi\in
\Hom(D_{j'}^{(h')},D_{j}^{(j)})^{*}\sub
\Hom(\TD_1,\TD_1)^{*},
$$
then for $1\leq h'\leq\min(j,j')$ we have
$$
\widetilde\chi(\tde_1\tga_1)=
\chi(\pi_{D_{j}^{(j)}}(\tde_1\tga_1)|_{D_{j'}^{(h')}})=
\chi(q_{h'\ra j}p_{j'\ra h'}),
$$
which are all the generators of the algebra $\RR(\La^c(v,d))^{G(V)}$
according to the Theorem \ref{invtheorem}.
\end{proof}

\lemm{} The map $\tphi: \fM(v,d) \to\TT_{x}^{a}$ is proper
and injective.

\begin{proof} We have the following diagrams
\begin{equation}\label{sdiagram}
\begin{CD}
\fM(v,d)@>{\tphi}>> \TT_x^a \\
@V{p}VV @V{\mmm_a}VV \\
\fM_0(v,d)@>{\phi_0}>> \fM_0(\tv,\td)
\end{CD}
\qquad
\begin{CD}
\fM(v,d)@>{\tphi}>> \TT_x^a \\
@V{p}VV @V{\mmm_a}VV \\
\fM_1(v,d)@>{\phi}>>  T_{x,\mu}
\end{CD}
\end{equation}
Since $\phi$
is a closed immersion and the morphisms $p$ and $\mmm_a$ are
projective, we see that $\tphi$ is proper.
Since all orbits in $\La^c_s(v,d)$
and $\La^c_s(\tv,\td)$ are closed, $\tphi$ is injective.
\end{proof}

\lemm{} The map $\tphi: \fM(v,d) \to\TT_x^a$ is an
isomorphism of algebraic varieties.

\begin{proof} Since $\tphi$ is a proper
injective morphism between connected smooth varieties of the same
dimension, $\tphi$ is an analytic isomorphism and therefore
an algebraic isomorphism.
\end{proof}

\slemm{} The map $\phi:\fM_1(v,d)\to T_{x,\mu}$ is an isomorphism
of algebraic varieties.

\begin{proof} Since $\mmm_a$ is surjective, from the diagram
(\ref{sdiagram}) we see that $\phi$ is surjective.
Since $\phi$ is a surjective closed immersion,
and both $\fM_1(v,d)$ and $T_{x,\mu}$ are reduced varieties
over $\C$, $\phi$ is an algebraic isomorphism.
\end{proof}

\sus{The immersions $\psi$ and $\tpsi$}
These immersions were constructed in section \ref{globalpsi}.

\se{Application to representation theory: $(\gl(n),\gl(m))$-duality}\label{applicationsRepTheory}
The relationship between quiver varieties and affine Grassmannians
provides a natural framework for $(GL(n),GL(m))$ duality.

\sus{Skew $(GL(n),GL(m))$ duality}

\sss{} Let $V=\C^m$ and $W=\C^n$ be two vector
spaces. Let us consider the $\gl(m)\times \gl(n)$ bimodule $V\otimes W$
and its $N$-th exterior power $\wedge^N(V\otimes W)$.
We have the following decomposition \cite[4.1.1]{H}:
\begin{equation}\label{howew}
\wedge^N(V\otimes W)=\bigoplus_{\la}V_\la\otimes W_{\cla},
\end{equation}
where $\la$ are all partitions of $N$ which fit into the $n\times m$
box, $V_\la$ is the highest weight representation of $\gl(m)$
with highest weight $\la$
and $W_{\cla}$ is the highest weight representation of $\gl(n)$
with highest weight $\cla$.

\sss{} Considering
$V\otimes W$ as a $\gl(m)$ module $V\otimes\C^n$, we have the following
decomposition:
\begin{equation}\label{wedge}
\wedge^N(V\otimes W)=\bigoplus_{a_1+\dots +a_n=N}
\wedge^{a_1}V\otimes\dots\otimes\wedge^{a_n}V.
\end{equation}

Considered as a representation of the torus
$(\C^{\times})^n\sub \gl(n)$ the vector space
$\wedge^{a_1}V\otimes\dots\otimes\wedge^{a_n}V$ has weight
$a=(a_1,\dots,a_n)$. Thus decompositions (\ref{howew}) and (\ref{wedge})
imply the following formula
\begin{equation}\label{weightmult1}
\Hom_{\gl(m)}
(\wedge^{a_1}V\otimes\dots\otimes \wedge^{a_n}V,
V_{\la}) \simeq W_{\cla}(a),
\end{equation}
where $W_{\cla}(a)$ is the weight space corresponding to weight $a$
of the $\gl(n)$ highest weight module $W_{\cla}$.

\sss{Geometric skew duality}
We construct
a based version of the isomorphism (\ref{weightmult1}),
i.e., a geometric skew $(GL(n),GL(m))$ duality.
More precisely, with $N,v,d,a,\la$ as in \ref{qdatatogl},
we identify
the right hand side
with
$\HH(\pi^{-1} (L_\la))$, where $L_{\la}$ is a lattice in
the affine Grassmannian $\GG$,
and the left hand side
with
$\HH(\fL(v,d))$ by Theorem \ref{theonakajima}.
The identification of irreducible components
$\Irr\pi^{-1} (L_\la)=\Irr\fL(v,d)$,
which follows from the isomorphism (\ref{fiberiso})
matches the
natural basis of the space
of intertwiners
$\Hom_{GL(m)}
(\wedge^{a_1}V\otimes\dots\otimes \wedge^{a_n}V,
V_{\la})$
arising from the affine Grassmannian construction
(i.e., $\Irr\pi^{-1} (L_\la)$), and the
natural basis of the weight space $W_{\cla}(a)$ in the Nakajima
construction
(i.e., $\Irr\fL(v,d)$). Altogether:
$$
\Hom_{GL(m)}
(\wedge^{a_1}V\otimes\dots\otimes \wedge^{a_n}V,V_{\la})
\simeq\HH(\pi^{-1} (L_\la))
\simeq\HH(\fL(v,d))
\simeq W_{\cla}(a).
$$

\sss{} Dually, we have
\begin{equation}\label{weightmult2}
\Hom_{\gl(n)}
(\wedge^{c_1}W \otimes\dots\otimes \wedge^{c_m}W, W_{\cla})=V_{\la}(c),
\end{equation}
where $V_{\la}(c)$ is the weight space corresponding to the
weight $c=(c_1\dots, c_m)$
of the $\gl(m)$ highest weight module $V_{\la}$.

\sus{Symmetric $(GL(m),GL(m))$ duality}

\sss{} Analogously, if we consider the $N$-th symmetric power
$\Sym^N(V\otimes V)$
of the $\gl(m)\times \gl(m)$ bimodule $V\otimes V$, we have the following
decomposition (a particular case of \cite[2.1.2]{H}):
\begin{equation}\label{howes}
\Sym^N(V\otimes V)=\bigoplus_{\la}V_\la\otimes V_{\la},
\end{equation}
where the sum is over all partitions $\la$ of $N$ with
at most $m$ parts.

Considering $V\otimes V$ as a $\gl(m)$ module $V\otimes\C^m$, we
have the following
decomposition:
\begin{equation}\label{sym}
\Sym^N(V\otimes V)=\bigoplus_{c_1+\dots +c_m=N}
\Sym^{c_1}V\otimes\dots\otimes\Sym^{c_m}V.
\end{equation}
Thus decompositions (\ref{howes}) and (\ref{sym})
imply the following formula
\begin{equation}\label{weightmultsym}
\Hom_{\gl(m)}
(\Sym^{c_1}V\otimes\dots\otimes \Sym^{c_m}V,
V_{\la})=V_{\la}(c),
\end{equation}
where $V_{\la}(c)$ is the weight space corresponding to weight $c$
of the $\gl(m)$ highest weight module $V_{\la}$.

\sss{} Combining the equations (\ref{weightmult2}) and
(\ref{weightmultsym}) we get
\begin{equation}\label{weightmultmix}
\Hom_{\gl(n)}
(\wedge^{c_1}W \otimes\dots\otimes \wedge^{c_m}W,
W_{\cla})=
\Hom_{\gl(m)}
(\Sym^{c_1}V\otimes\dots\otimes \Sym^{c_m}V,
V_{\la}).
\end{equation}

\sss{Geometric symmetric duality} Geometry allows us to find a \emph{based}
isomorphism of the left and right hand side of (\ref{weightmultmix}).
Let $N,v,d,a,\la$ be as in \ref{qdatatogl}. First of all it follows from the quiver tensor product
constructions of Malkin \cite{Mal} and Nakajima \cite{N01b}
that the relevant irreducible components
$\Irr \tfg^{n,c}_x$ of the Spaltenstein fiber
over a nilpotent of type $\la$
index a natural basis in the left hand side of (\ref{weightmultmix}).
Here
$$
\tfg^{n,c}_x=\{(x,F)\in \\gl(D)\times \FF^{m,c}~|~x(F_i)\sub F_i
\text{ and } x \text{ acts on } F_i/F_{i-1}
\text{ as a regular nilpotent }\}.
$$
Now consider another convolution Grassmannian:
\begin{equation}
\begin{split}
\TGG^c & =\barr\GG_{c_1\om_1}\ast\cdots\ast\barr\GG_{c_m\om_1} \\
& =
\{\ L_0\sub\ L_1\sub\cddd\sub L_n ~
| ~
\dim L_i/L_{i-1}=c_i,\  z|_{L_i/L_{i-1}}
\text{ is a regular nilpotent }\},
\end{split}
\end{equation}
where $\om_1$ is the first fundamental weight of $GL(m)$.
We have a map $\pi: \TGG^c\to \GG$ defined by
$\pi: (L_0\sub L_1\sub \dots \sub L_n)\mapsto L=L_n$.
Consider $\pi^{-1}(L_\la)$ for $L_\la\in\GG$.
It follows from the Geometric Satake Correspondence that
the set of relevant irreducible components $\Irr\pi^{-1} (L_\la)$
indexes a basis in the right hand side of (\ref{weightmultmix}).

It is clear that the varieties
$\tfg^{n,c}_x\simeq \pi^{-1}(L_\la)$ are isomorphic.
This isomorphism gives us a bijection
$\Irr \tfg^{n,c}_x=\Irr\pi^{-1} (L_\la)$.

Summarizing:
\begin{equation}\nonumber
\begin{split}
\Hom_{GL(n)}
(\wedge^{c_1}W\otimes\dots\otimes \wedge^{c_m}W,
W_{\cla}) & \simeq \HH(\tfg^{n,c}_x) \\
& \simeq\HH(\pi^{-1} (L_\la)) \\
& \simeq
\Hom_{GL(m)}
(\Sym^{c_1}V\otimes\dots\otimes \Sym^{c_m}V,
V_{\la}).
\end{split}
\end{equation}

\sss{Remark} The second author has greatly benefited from a class taught by W.~ Wang at Yale \cite{Wa1}.
The ``geometric symmetric duality'' above has a lot in common with the construction described in \cite{Wa2} and
we believe that the ``geometric skew duality'' construction answers a question posed by Weiqiang Wang.

\end{document}